\documentclass[final,leqno]{siamltex}
\usepackage{bbm,xcolor,bm}
\usepackage{amssymb}
\usepackage{amsfonts}
\usepackage{mathrsfs}
\usepackage{multirow}
\usepackage{graphicx}
\usepackage{rotating}
\usepackage{subfigure}
 \usepackage{color}
\usepackage{cite}
\usepackage[cmex10]{amsmath}
\usepackage{epstopdf, url, ulem}

\newcommand{\Lone}{\textnormal{L}_1}
\newcommand{\Ltwo}{\textnormal{L}_2}
 \usepackage{algorithm}

\newcommand{\R}{\mathbf{R}}
\newcommand{\B}{\mathbf{B}}

\newcommand{\x}{\mathbf{x}}
\newcommand{\y}{\mathbf{y}}

\newcommand{\I}{\mathbf{I}}

\newcommand{\J}{\mathcal{J}}

\newcommand{\nvar}{\textsc{n}_\textnormal{var}}

\def\x{\mathbf{x}}
\def\u{\mathbf{u}}
\def\y{\mathbf{y}}
\def\z{\mathbf{z}}
\def\R{\mathbf{R}}
\def\B{\mathbf{B}}
\def\HH{\mathbf{H}}

\def\nens{\textsc{n}_{\rm ens}}

\usepackage{algpseudocode,caption}
\algblock{ParFor}{EndParFor}
\algnewcommand\algorithmicparfor{\textbf{For all}}
\algnewcommand\algorithmicpardo{\textbf{do in parallel}}
\algnewcommand\algorithmicendparfor{\textbf{end\ For all}}
\algrenewtext{ParFor}[1]{\algorithmicparfor\ #1\ \algorithmicpardo}
\algrenewtext{EndParFor}{\algorithmicendparfor}

\newcommand{\pk}{{\{k\}}}
\newcommand{\pkp}{{\{k+1\}}}

\newcommand{\myfigurewidth}{0.475\textwidth}
\newcommand{\myfigureheight}{0.35\textwidth}
\title{Robust data assimilation using $\Lone$ and Huber norms\thanks{This work was supported by AFOSR DDDAS program through the award AFOSR FA9550--12--1--0293--DEF managed by Dr. Frederica Darema, and by the Computational Science laboratory at Virginia Tech. M. Ng's research is supported in part by HKRGC GRFs 202013 and 12301214, and HKBU FRG2/14-15/087.}}

\author{Vishwas Rao \footnotemark[2] \footnotemark[4]
\and Adrian Sandu \footnotemark[3] \footnotemark[4]
\and Michael Ng \footnotemark[5] \footnotemark[6] \and \newline
Elias D. Nino-Ruiz \footnotemark[7] \hspace{0.1mm} \footnotemark[4] }

\begin{document}
\thispagestyle{empty}
\setcounter{page}{0}

\makeatletter
\def\Year#1{%
  \def\yy@##1##2##3##4;{##3##4}%
  \expandafter\yy@#1;
}
\makeatother

\begin{Huge}
\begin{center}
Computational Science Laboratory Technical Report CSL-TR-\Year{\the\year}-{\tt 21} \\
\today
\end{center}
\end{Huge}
\vfil
\begin{huge}
\begin{center}
{\tt Vishwas Rao, Adrian Sandu, Michael Ng and Elias Nino-Ruiz}
\end{center}
\end{huge}

\vfil
\begin{huge}
\begin{it}
\begin{center}
``{\tt Robust data assimilation using $\Lone$ and Huber norms}''
\end{center}
\end{it}
\end{huge}
\vfil

\begin{large}
\begin{center}
Computational Science Laboratory \\
Computer Science Department \\
Virginia Polytechnic Institute and State University \\
Blacksburg, VA 24060 \\
Phone: (540)-231-2193 \\
Fax: (540)-231-6075 \\ 
Email: \url{sandu@cs.vt.edu} \\
Web: \url{http://csl.cs.vt.edu}
\end{center}
\end{large}

\vspace*{1cm}

\begin{tabular}{ccc}
\includegraphics[width=2.5in]{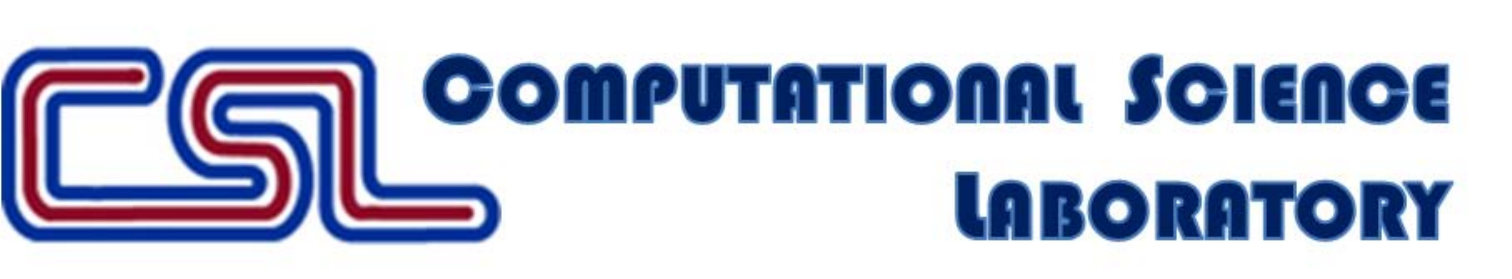}
&\hspace{2.5in}&
\includegraphics[width=2.5in]{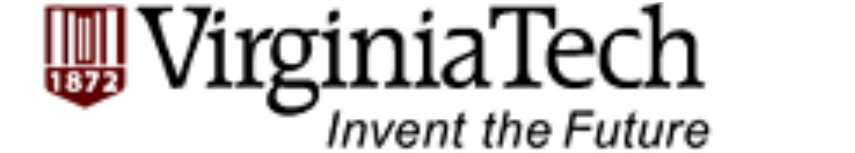} \\
{\bf\em\large Compute the Future} &&\\
\end{tabular}

\newpage

\maketitle
\renewcommand{\thefootnote}{\fnsymbol{footnote}}
\footnotetext[2]{E-mail: visrao@vt.edu}
\footnotetext[3]{E-mail: sandu@cs.vt.edu}
\footnotetext[4]{Computational Science Laboratory, Department of Computer Science, Virginia Tech.  2202 Kraft Drive, Blacksburg, Virginia 24060.}
\footnotetext[5]{E-mail: mng@math.hkbu.edu.hk}
\footnotetext[6]{Department of Mathematics, Hong Kong Baptist University}
\footnotetext[7]{E-mail: enino@vt.edu}
\renewcommand{\thefootnote}{\arabic{footnote}}

\begin{abstract}
Data assimilation is the process to fuse information from priors, observations of nature, and  numerical models, in order to obtain best estimates of the parameters or state of a physical system of interest. Presence of large errors in some observational data, e.g., data collected from a faulty instrument, negatively affect the quality of the overall assimilation results.
This work develops a systematic framework for robust data assimilation. The new algorithms continue to produce good analyses in the presence of observation outliers. The approach is based on replacing the traditional $\Ltwo$ norm formulation of data assimilation problems with formulations based on $\Lone$ and Huber norms. Numerical experiments using the Lorenz-96 and the shallow water on the sphere models illustrate how the new algorithms outperform traditional data assimilation approaches in the presence of data outliers.
\end{abstract}

\section{Introduction} \label{sec:intro}
Dynamic data-driven application systems  (DDDAS \cite{dddas}) integrate computational simulations and physical measurements in symbiotic and dynamic feedback control systems. Within the DDDAS paradigm, data assimilation (DA) defines a class of inverse problems that fuses information from an imperfect computational model based on differential equations (which encapsulates our knowledge of the physical laws that govern the evolution of the real system), from noisy observations (sparse snapshots of  reality), and from an uncertain prior (which encapsulates our current knowledge of reality). Data assimilation integrates these three sources of information and the associated uncertainties in a Bayesian framework to provide the posterior, i.e., the probability distribution conditioned on the uncertainties in the model and observations.

Two approaches to data assimilation have gained widespread popularity: ensemble-based estimation and variational methods. The ensemble-based methods are rooted in statistical theory, whereas the variational approach is derived from optimal control theory. The variational approach formulates data assimilation as a nonlinear optimization problem constrained by a numerical model. The initial conditions (as well as boundary conditions, forcing, or model parameters) are adjusted to minimize the discrepancy between the model trajectory and a set of time-distributed observations. In real-time operational settings the data assimilation process is performed in cycles: observations within an assimilation window are used to obtain an optimal trajectory, which provides the initial condition for the next time window, and the process is repeated in the subsequent cycles.

Large errors in some observations can adversely impact the overall solution to the data assimilation system, e.g., can lead to spurious features in the analysis \cite{lorenc1986analysis}. Various factors contribute to uncertainties in observations. Faulty and malfunctioning sensors are a major source of large uncertainties in observations, and data quality control (QC) is an important component of any weather forecasting  DA system  \cite{lorenc1988objective}. The goal of QC is to ensure that only correct observations are used in the DA system. Erroneous observations can lead to spurious features in the resulting solution \cite{lorenc1986analysis}. Departure of the observation from the background forecast is usually used as a measure for accepting or rejecting observations: observations with large deviations from the background forecast are rejected in the quality control step \cite{hollingsworth1986monitoring}. However, this process has important drawbacks. Tavolato and Isakksen present several case studies \cite{tavolato2014use} that demonstrate that rejecting observations on the basis of background departure statistics leads to the inability to capture small scales  in the analysis.

As an alternative to rejecting the observations which fail the QC tests, Tavolato and Isakksen \cite{tavolato2014use} propose the use of Huber norm in the analysis. Their Huber norm approach adjusts error covariances of observations based on their background departure statistics. No observations are discarded; but smaller weights are given to low quality observations  (that show large departures from the background). Huber norm has been used in the context of robust methods for seismic inversion \cite{Guitton_2003_robust-inversion} and ensemble Kalman filtering \cite{roh2013observation}. The issue of outliers has also been tackled by regularization techniques using $\Lone$ norms \cite{ebtehaj2013variational}.
To the best of our knowledge, no prior work has fully addressed the construction of robust data assimilation methods using these norms. For example \cite{tavolato2014use} adjusts the observation data covariances but then applies traditional assimilation algorithms.

This paper develops a systematic mathematical framework for robust data assimilation. Different data assimilation algorithms are reformulated as optimization problems using $\Lone$ and Huber norms, and are solved via the alternating direction method of multipliers (ADMM) developed by Boyd \cite{Boyd_2011_ADMM}, and via the half-quadratic minimization method presented in \cite{nikolova2005analysis}.

The remainder of the paper is organized as follows: Section \ref{sec:DA} reviews the data assimilation problem and the traditional solution algorithms. Section \ref{sec:3DVar_robust} presents the new robust algorithms for 3D-Var, Section \ref{sec:Robust-4DVar} for 4D-Var, and Section \ref{sec:L1-EnKF} develops new robust algorithms for ensemble based approaches. Numerical results are presented in Section \ref{sec:numexp}. Concluding remarks and future directions are discussed in Section \ref{sec:conc}.

\section{Data assimilation}
\label{sec:DA}
Data assimilation (DA) is the process of fusing information from priors, imperfect model predictions, and noisy data, to obtain a consistent description of the true state $\x^{\rm true}$ of a physical system \cite{daley1993, kalnay2003, sandu2011chemical, sandu2005adjoint}. The resulting best estimate is called the analysis $\x^{\rm a}$.

The prior information encapsulates our current knowledge of the system. Prior information typically consists of a background estimate of the state $\x^{\rm b}$, and the corresponding background error covariance matrix $\mathbf{B}$.

The model captures our knowledge about the physical laws that govern the evolution of the system. The model evolves an initial state $\x_0 \in \mathbb{R}^n$ at the initial time $t_{ 0 }$ to states $\x_{i} \in \mathbb{R}^n$ at future times $t_{i}$. A general model equation is:
\begin{equation}
\label{eqn:genmodel}
 \x_{i} = \mathcal{M}_{t_0 \rightarrow t_{i}} \left(\x_0\right), \quad  i=1,\cdots, N.
\end{equation}
Observations are noisy snapshots of reality available at discrete time instances. Specifically, measurements $\y_{i} \in \mathbb{R}^m$ of the physical state $\x^{\rm true}\left(t_{i}\right)$ are taken at times $t_{i}$,
\begin{equation}
\label{eqn:genobservation}
\y_{i} = \mathcal{H}(\x_i) + \varepsilon_i, \quad  \varepsilon_i \sim \mathcal{N}(\mathbf{0},\R_i),
\quad i=1,\cdots, N,
\end{equation}
where the observation operator $\mathcal{H} : \mathbb{R}^n \to \mathbb{R}^m$ maps the model state space onto the observation space. The random observation errors $\varepsilon_i$ are assumed to have normal distributions.

Variational methods solve the DA problem in an optimal control framework. Specifically, one adjusts a control variable (e.g., model parameters) in order to minimize the discrepancy between model forecasts and observations. Ensemble based methods are sequential in nature. Generally, the error statistics of the model estimate is represented by an ensemble of model states and ensembles are propagated to predict the error statistics forward in time. When the observations are assimilated, the analysis is obtained by operating directly on the model states. It is important to note that all the approaches discussed in the following subsections aim to minimize the $\Ltwo$ norm of the discrepancy between observations and model predictions.

Ensemble-based methods employ multiple model realizations in order to approximate the background error distribution of any model state. The empirical moments of the ensemble are used in order to estimate the background state and the background error covariance matrix. Since the model space is several times larger than the ensemble size, localization methods are commonly used in order to mitigate the impact of sampling errors (spurious correlations). For example, the local ensemble transform Kalman filter \cite{Hunt_2007_LETKF} performs the assimilation of observations locally in order to avoid the impact of analysis innovations coming from distant components.

We next review these traditional approaches.

\subsection{Three-dimensional variational data assimilation}
The three dimensional variational (3D-Var) DA approach processes the observations successively at times $t_{1},t_{2}, \dots, t_{N}$. The background state (i.e., the prior state estimate at time $t_{i}$) is given by the model forecast started from the previous analysis (i.e., the best state estimate at $t_{i-1}$):
\begin{equation}\label{eqn:3D-VarBg}
 \x_{i}^{\rm b} = \mathcal{M}_{t_{i-1} \rightarrow t_{ i }} \left(\x^{\rm a}_{i-1}\right)\,.
\end{equation}
The discrepancy between the model state $\x_i$ and observations at time $t_{i}$, together with the departure of the state from the model forecast $\x_{i}^{\rm b}$, are measured by the 3D-Var cost function
\begin{subequations}
\label{eqn:3DVar}
\begin{equation}
\label{eqn:3DVarCF}
 \mathcal{J}\left(\x_i\right) = \frac{1}{2}\| \x_i - \x_i^{\rm b} \|_{\B_i^{-1}}^2 +
\frac{1}{2} \, \| \mathcal{H}(\x_i) - \y_i \|_{\R_i^{-1}}^2.
\end{equation}
When both the background and observation errors are Gaussian the function \eqref{eqn:3DVarCF} equals the negative logarithm of the posterior probability density function given by Bayes' theorem. While in principle a different background covariance matrix should be used at each assimilation time, in practice the same matrix is re-used throughout the assimilation window, $\B_i =\B, \,i=1,\dots,N$. The 3D-Var analysis is the maximum a-posteriori (MAP) estimator, and is computed as the state that minimizes \eqref{eqn:3DVarCF}
\begin{equation}
\label{eqn:3DVarL2}
 \x_i^{\rm a} = \arg \min_{\x_i}\, \mathcal{J}\left(\x_i\right).
\end{equation}
\end{subequations}

\subsection{Four-dimensional variational data assimilation}
Strong-constraint four-dimensional variational (4D-Var) DA processes simultaneously all observations at all times $t_{1},t_{2}, \dots, t_{N}$ within the assimilation window. The control parameters are typically the initial conditions $\x_0$, which uniquely determine the state of the system at all future times under the assumption that the model \eqref{eqn:genmodel} perfectly represents reality. The background state is the prior best estimate of the initial conditions $\x_0^{\rm b}$, and has an associated initial background error covariance matrix $\mathbf{B}_0$. The 4D-Var problem provides the estimate $\x_0^{\rm a}$ of the true initial conditions as the solution of the following optimization problem
\begin{subequations}
\label{eqn:L2-4dvar}
\begin{eqnarray}
\label{eqn:ip}
&&~ \x_0^{\rm a}  =  \underset{\x_0} {\text{ arg\, min}}~~ \J\left(\x_0\right) \qquad
 \text{subject to}~ \text{\eqref{eqn:genmodel}},\\
\label{eqn:fdvar-cf}
&&~ \mathcal{J}\left(\x_0\right) = \frac{1}{2}\| \x_0 - \x_0^{\rm b} \|_{\B_0^{-1}}^2 +
\frac{1}{2} \, \sum_{i=1}^{N}\| \mathcal{H}(\x_i) - \y_i \|_{\R_i^{-1}}^2 \,.
\end{eqnarray}
\end{subequations}
The first term of the sum \eqref{eqn:fdvar-cf} quantifies the departure of the solution $\x_0$ from the background state $\x_0^{\rm b}$ at the initial time $t_0$. The second term measures the mismatch between the forecast trajectory (model solutions $\x_{i}$) and observations $\y_{i}$ at all times $t_{i}$ in the assimilation window. The covariance matrices $\mathbf{B}_0$ and $\mathbf{R}_{i}$ need to be predefined, and their quality influences the accuracy of the resulting analysis.
Weak-constraint 4D-Var \cite{sandu2011chemical} removes the perfect model assumption by allowing a model error
$\eta_{i+1} = \x_{i+1} - \mathcal{M}_{t_i \to t_{i+1}} \left(\x_{i}\right)$.

In this paper we focus on the strong-constraint 4D-Var formulation \eqref{eqn:ip}. The minimizer of \eqref{eqn:ip} is computed iteratively using gradient-based numerical optimization methods. First-order adjoint models provide the gradient of the cost function \cite{cacuci2005sensitivity}, while second-order adjoint models provide the Hessian-vector product (e.g., for Newton-type methods) \cite{Sandu_2008_SOA}. The methodology for building and using various adjoint models for optimization, sensitivity analysis, and uncertainty quantification is discussed in \cite{sandu2005adjoint, cioaca2012second}. Various strategies to improve the the 4D-Var data assimilation system are described in \cite{cioaca2014optimization}. The procedure to estimate the impact of observation and model errors is developed in \cite{rao2014posterioriJournal, rao2014posteriori}. A framework to perform derivative free variational data assimilation using the trust-region framework is given in \cite{ruiz2015derivative}.

\subsection{The ensemble Kalman filter}
The uncertain background state at time $t_i$ is assumed to be normally distributed $\x_i \sim \mathcal{N}\left( \overline{\x}_i^\textnormal{b},\mathbf{P}^{\rm b}_i \right)$. Given the new observations $\y_i$ at time $t_i$, the Kalman filter computes the Kalman gain matrix, the analysis state (i.e., the posterior mean), and the analysis covariance matrix (i.e., the posterior covariance matrix) using the following formulas, respectively:
\begin{subequations}
\label{eqn:kf}
\begin{eqnarray}
\label{eqn:kf-gain}
\mathbf{K}_i &=& \mathbf{P}^{\rm b}_i\, \HH_i^\mathrm{T} \left( \HH_i\, \mathbf{P}^{\rm b}_i \, \HH_i^\mathrm{T} + \R_i \right)^{-1}, \\
\label{eqn:kf-mean}
\overline{\x}_i^\textnormal{a} &=& \overline{\x}_i^\textnormal{b} + \mathbf{K}_i\, \left(  \y_i - \mathcal{H}(\overline{\x}_i^\textnormal{b}) \right), \\
\label{eqn:kf-analysis}
\mathbf{P}^{\rm a}_i &=& \left( \I - \mathbf{K}_i \, \HH_i \right)\, \mathbf{P}^{\rm b}_i.
\end{eqnarray}
\end{subequations}

The ensemble Kalman filter describes the background probability density at time $t_i$ by an ensemble of $\nens$ forecast states
$
\{ \x_i^{ {\rm b}\langle\ell\rangle}\}_{\ell = 1,\dots,\nens}.
$
The background mean and covariance are estimated from the ensemble:
\begin{equation}
\label{eqn:enkf-mean-cov}
\overline{\x}^{\rm b}_i \approx \frac{1}{\nens} \sum_{\ell=1}^{\nens} \x_i^{ {\rm b}\langle\ell\rangle},\qquad
\quad \mathbf{P}^{\rm b}_i \approx \frac{1}{\nens-1} \, \mathbf{X}^{\rm b}_i\, (\mathbf{X}^{\rm b}_i)^T,
\end{equation}
where the matrix of state deviations from the mean is
\[
\mathbf{X}^{\rm b}_i = \left[\x_i^{ {\rm b}\{1\}}-\overline{\x}^{\rm b}_i, \dots, \x_i^{ {\rm b}\{\nens\}} -\overline{\x}^{\rm b}_i \right]
\approx  \frac{1}{\sqrt{\nens-1}}  \left( \mathbf{P}^{\rm b}_i \right)^{1/2}.
\]
The background observations and their deviations from the mean are defined as:
\[
\overline{\mathbf{h}}^{\rm b}_i = \frac{1}{\nens} \sum_{\ell=1}^{\nens} \mathcal{H}\left( \x_i^{ {\rm b}\langle\ell\rangle} \right),\qquad
\mathbf{Y}^{\rm b}_i = \left[ \mathcal{H}\left( \x_i^{ {\rm b}\{1\}} \right)-\overline{\mathbf{h}}^{\rm b}_i , \dots, \mathcal{H}\left( \x_i^{ {\rm b}\{\nens\}} \right) - \overline{\mathbf{h}}^{\rm b}_i \right].
\]
Any vector in the subspace spanned by the ensemble members can be represented as
\[
\x = \overline{\x}^{\rm b}_i + \mathbf{X}^{\rm b}_i \, w.
\]

The ensemble Kalman filter makes the approximations \eqref{eqn:enkf-mean-cov}
%
and computes the Kalman gain matrix \eqref{eqn:kf-gain} as follows:
\begin{subequations}
\begin{eqnarray}
\mathbf{K}_i
&=& \mathbf{X}^{\rm b}_i\,\mathbf{S}_i \, (\mathbf{Y}^{\rm b}_i)^\mathrm{T}\, \R_i^{-1}, \\
\label{eqn:letkf-special-matrix}
\mathbf{S}_i &:=& \left( (\nens-1) \mathbf{I} +  \left(\mathbf{Y}^{\rm b}_i\right)^T\, \R_i^{-1}\, \mathbf{Y}^{\rm b}_i \right)^{-1}.
\end{eqnarray}
\end{subequations}
The mean analysis \eqref{eqn:kf-mean} is:
\begin{eqnarray}
\label{eqn:enkf-mean}
\overline{\x}_i^\textnormal{a}
&=& \overline{\x}_i^\textnormal{b} + \mathbf{X}_i^\textnormal{b}\,\mathbf{S}_i \, (\mathbf{Y}_i^\textnormal{b})^\mathrm{T}\, \R_i^{-1} \, \left(  \y_i - \mathcal{H}(\overline{\x}_i^\textnormal{b}) \right)
 \end{eqnarray}
Using $\overline{\x}_i^\textnormal{a} = \overline{\x}_i^\textnormal{b} + \mathbf{X}_i^\textnormal{b} \, \overline{w}^\textnormal{a}_i$ we have the mean analysis in ensemble space:
\begin{eqnarray}
\label{eqn:enkf-analysis-mean}
\overline{w}^\textnormal{a}_i
&=& \mathbf{S}_i \, (\mathbf{Y}_i^\textnormal{b})^\mathrm{T}\, \R_i^{-1} \, \left(  \y_i - \mathcal{H}(\overline{\x}_i^\textnormal{b}) \right).
\end{eqnarray}
The analysis error covariance \eqref{eqn:kf-analysis} becomes:
\begin{eqnarray}
\label{eqn:enkf-analysis-covariance}
\mathbf{P}^{\rm a}_i
&=& \mathbf{X}_i^\textnormal{b} \,
\mathbf{S}_i \, (\mathbf{X}_i^\textnormal{b})^T.
\end{eqnarray}

\paragraph{Ensemble square root filter (EnSRF)}

The 3D-Var cost function \eqref{eqn:3DVarCF} is formulated with the background error covariance  $\B_i$ replaced by the forecast ensemble covariance matrix $\mathbf{P}^{\rm b}_i$ \eqref{eqn:enkf-mean-cov}.
It can be shown that the analysis mean weights \eqref{eqn:enkf-analysis-mean} are the minimizer of \cite{Hunt_2007_LETKF}:
\begin{equation}
\label{eqn:EnKF-L2}
\overline{w}^\textnormal{a}_i := \arg\min_w ~\mathcal{J}(w), \quad \mathcal{J}(w)= (\nens-1)\, \| w \|_2^2 +
\| \mathcal{H}\left(\overline{\x}^{\rm b}_i + \mathbf{X}^{\rm b}_i \, w\right) - \y_i \|_{\R_i^{-1}}^2.
\end{equation}
%

The local transform ensemble Kalman filter \cite{Hunt_2007_LETKF} computes the symmetric square root of matrix \eqref{eqn:letkf-special-matrix}
\begin{equation}
\label{eqn:letkf-special-matrix-factorization}
\mathbf{S}_i  = (\nens-1)^{-1}\, \mathbf{W}_i\, \mathbf{W}_i^T,
\end{equation}
and obtains the analysis ensemble weights by adding columns of the factors to the mean analysis weight:
\begin{equation}
\label{eqn:letkf-analysis-ensemble}
w_i^{\textnormal{a}\langle\ell\rangle} = \overline{w}^\textnormal{a}_i + \mathbf{W}_i(:,\ell),
\qquad
\x_i^{\textnormal{a}\langle\ell\rangle} = \overline{\x}_i^\textnormal{b} + \mathbf{X}_i^\textnormal{b} \, w_i^{\textnormal{a}\langle\ell\rangle}.
\end{equation}

\paragraph{Traditional ensemble Kalman filter (EnKF)}
In the traditional version of EnKF the filter \eqref{eqn:enkf-mean} is applied to each ensemble member:
\begin{eqnarray*}
\x_i^{\textnormal{a}\langle\ell\rangle}
&=&  \x_i^{\textnormal{b}\langle\ell\rangle} + \mathbf{X}_i^\textnormal{b}\,\mathbf{S}_i \, (\mathbf{Y}_i^\textnormal{b})^\mathrm{T}\, \R_i^{-1} \, \left(  \y_i^{\langle\ell\rangle} - \mathcal{H}(\x_i^{\textnormal{b}\langle\ell\rangle}) \right),
\quad \ell = 1,\dots,\nens,
 \end{eqnarray*}
where $\y_i^{\langle\ell\rangle}$ are perturbed observations. This leads to the analysis weights
\begin{equation}
\label{eqn:enkf-analysis-member}
\begin{split}
w_i^{\textnormal{a}\langle\ell\rangle}
=&~  w_i^{\textnormal{b}\langle\ell\rangle} + \mathbf{S}_i \, (\mathbf{Y}_i^\textnormal{b})^\mathrm{T}\, \R_i^{-1} \, \left(  \y_i^{\langle\ell\rangle} - \mathcal{H}(\x_i^{\textnormal{b}\langle\ell\rangle}) \right) \\
\approx&~  \mathbf{S}_i \,  \left(
(\nens-1)\,w_i^{\textnormal{b}\langle\ell\rangle} + (\mathbf{Y}_i^\textnormal{b})^\mathrm{T}\, \R_i^{-1} \, \left(  \y_i^{\langle\ell\rangle} - \mathcal{H}(\overline{\x}_i^\textnormal{b})  \right)
\right),
\end{split}
\end{equation}
where the second relation comes from linearizing the observation operator about the background mean state.
%
%
A comparison of \eqref{eqn:enkf-analysis-member} with \eqref{eqn:enkf-analysis-mean} reveals that the classical EnKF solution \eqref{eqn:enkf-analysis-member} obtains the weights of the $\ell$-th ensemble member by minimizing the cost function:
\begin{equation}
\label{eqn:enkf-traditional-opt-L2}
\begin{split}
w_i^{\textnormal{a}\langle\ell\rangle} &= \arg\min_w ~ \mathcal{J}^{\langle\ell\rangle}(w), \\
\mathcal{J}^{\langle\ell\rangle}(w) &= (\nens-1)\, \| w - w_i^{\textnormal{b}\langle\ell\rangle} \|_2^2 +
\| \mathcal{H}\left(\overline{\x}^{\rm b}_i + \mathbf{X}^{\rm b}_i \, w\right) - \y^{\langle\ell\rangle}_i \|_{\R_i^{-1}}^2.
\end{split}
\end{equation}

The paper develops next a systematic framework to perform robust DA using Huber and $\Lone$ norms. This will ensure that important information coming from outliers is not rejected but used during the quality control.

\section{Robust 3D-Var data assimilation} \label{sec:3DVar_robust}

\subsection{$\Lone$-norm 3D-Var} \label{sec:L13DVar}
%
Consider the 3D-Var cost function in equation \eqref{eqn:3DVarCF} that penalizes the discrepancy between the model state $\x_i$ and observations at time $t_i$, together with the departure of the state from the model forecast $\x_i^{\rm b}$. Function \eqref{eqn:3DVarCF} represents the negative log-likelihood posterior probability density under the assumption that both the background and the observation errors are Gaussian. In particular, the scaled innovation vector
\begin{equation}
\label{eqn:define-z}
\z_i := \R_i^{-1/2}\, \left[ \mathcal{H}(\x_i) - \y_i \right],
\end{equation}
which represents scaled observation errors, is assumed to have a normal distribution $\z_i \sim \mathcal{N}(\mathbf{0},\mathbf{I})$ (each observation error $(\z_i)_\ell$ is independent and has a standard normal distribution). We are interested in the case where some observation errors have large deviations. To model this we assume that the scaled observation errors  \eqref{eqn:define-z} are independent, and each of them has a univariate Laplace distribution:
\[
\mathcal{P}(z) = (2\,\lambda)^{-1}\, \exp\left(\, - \lambda^{-1}\, \left| z - \mu \right|\, \right), \quad
{\tt E}\left[ z \right] = 0, \quad {\tt Var}\left[ z \right] = 2\,\lambda^2.
\]
The univariate Laplace distribution models an exponential decay on both sides away from the mean $\mu$. Assuming that observation errors are unbiased ($\mu=0$) and independent the probability distribution of the scaled innovation vector $\z_i$ is
\begin{eqnarray}
\label{eqn:laplace-product-distribution}
\mathcal{P}(\z_i) \propto 
\exp\left( - \lambda^{-1}\, \left\| \z_i \right\|_1 \right) , \quad
\tt{E}\left[ \z_i \right] = \mathbf{0}, \quad \tt{Cov}\left[ \z_i \right] = 2\,\lambda^2\, \mathbf{I}.
\end{eqnarray}
Under these assumptions the negative log-likelihood posterior function \eqref{eqn:3DVarCF} measures  the discrepancy of the model state with observations in the $\Lone$ norm. This leads to the following revised cost function:
\begin{equation}
\label{eqn:3dvar-L1}
\mathcal{J}(\x_i) = \frac{1}{2}
\| \x_i - \x_i^{\rm b} \|_{\B_i^{-1}}^2 +
\frac{1}{\lambda}
\left\| \R_i^{-1/2}\, \left[ \mathcal{H}(\x_i) - \y_i \right] \right\|_{1}\,.
\end{equation}
In this paper we choose $\lambda=2$, which leads to a cost function \eqref{eqn:3dvar-L1} that is similar to \eqref{eqn:3DVarCF}. This translates in an assumed variance of each observation error component equal to eight. This is no restriction of generality: the solution process discussed here can be applied to any value of $\lambda$. For example, one can choose $\lambda=1/\sqrt{2}$ for an error component variance of one.

Following the Alternating Direction Method of Multipliers (ADMM) \cite{Boyd_2011_ADMM} we use the scaled innovation vector \eqref{eqn:define-z} to obtain the $\Lone$-3D-Var problem:
\begin{equation}
\label{eqn:3dvar-L1-ADMM}
\begin{split}
& \min\, \mathcal{J}(\x_i,\z_i) = \frac{1}{2}
\| \x_i - \x_i^{\rm b} \|_{\B_i^{-1}}^2 +
\frac{1}{2} \| \z_i \|_{1} \\
& \textnormal{subject to} \quad \z_i = \R_i^{-1/2}\, [ \mathcal{H}(\x_i) - \y_i ].
\end{split}
\end{equation}
The augmented Lagrangian equation for \eqref{eqn:3dvar-L1-ADMM} is given by
\begin{eqnarray}
\label{eqn:3dvar-augmented-Lagrangian}
\mathcal{L} &=& \frac{1}{2} \| \x_i - \x_i^{\rm b} \|_{\B_i^{-1}}^2 +
\frac{1}{2}\, \| \z_i \|_{1} - \bm{\lambda}_i^T\, \left[ \R_i^{-1/2} [ \mathcal{H}(\x_i) - \y_i ] - \z_i \right] \\
\nonumber
&& + \frac{\mu}{2}\, \| \R_i^{-1/2} [ \mathcal{H}(\x_i) - \y_i ] - \z_i \|_2^2,
\end{eqnarray}
where $\bm{\lambda}_i$ is the Lagrange multiplier and $\mu$ is the penalty parameter (a positive number).
Using the identity
\begin{eqnarray}
\nonumber
&&  \frac{\mu}{2}\,\left\Vert \R_i^{-1/2} [ \mathcal{H}(\x_i) - \y_i ] - \z_i - \frac{\bm{\lambda}_i}{\mu} \right\Vert_2^2 =
 \frac{\mu}{2}\,\left\Vert \R_i^{-1/2} [ \mathcal{H}(\x_i) - \y_i ] - \z_i  \right\Vert_2^2 \\
 \label{eqn:z-identity}
 && \qquad\qquad + \frac{1}{2\mu}\, \left\Vert  \bm{\lambda}_i \right\Vert_2^2
 - \bm{\lambda}_i^T\, \left( \R_i^{-1/2} [ \mathcal{H}(\x_i) - \y_i ] - \z_i  \right),
\end{eqnarray}
the augmented Lagrangian \eqref{eqn:3dvar-augmented-Lagrangian} can be written as
\begin{equation}
\label{eqn:3dvar-augmented-Lagrangian-alt}
\qquad~~
\mathcal{L} = \frac{1}{2} \| \x_i - \x_i^{\rm b} \|_{\B_i^{-1}}^2 +
\frac{1}{2}\, \| \z_i \|_{1} - \frac{1}{2\mu}\, \left\Vert  \bm{\lambda} \right\Vert_2^2  + \frac{\mu}{2}\,\left\Vert \R_i^{-1/2} [ \mathcal{H}(\x_i) - \y_i ] - \z_i - \frac{\bm{\lambda}_i}{\mu} \right\Vert_2^2.
\end{equation}
The cost function \eqref{eqn:3dvar-augmented-Lagrangian-alt} can be iteratively minimized as follows:
\begin{enumerate}
\item[Step 1:] Initialize $\x_i^{\{0\}} = \x_i^{\rm b}\,$, $\z_i^{\{0\}} = \R_i^{-1/2}\, [ \mathcal{H}(\x_i^{\rm b}) - \y_i ]\,$, $\bm{\lambda}_i^{\{0\}} = \bm{0}\,$, $\mu^{\{0\}} = 1$
 \item[Step 2:] Start with $\x_i^{\{0\}}$, $\z_i^{\{0\}}$, $\bm{\lambda}_i^{\{0\}}$, and $\mu^{\{0\}}$.
 \item[Step 3:] Fix $\z_i^{\{k\}}$, $\bm{\lambda}_i^{\{k\}}$, and $\mu^{\{k\}}$, and solve
\begin{subequations}
\label{eqn:L1_modified_3dvar}
\begin{eqnarray}
\x^{\{k+1\}}_i &:=&
\arg\min_{\x}~
\frac{1}{2}
\left\Vert  \x - \x_i^{\rm b} \right\Vert_{\B_i^{-1}}^2 \\
\nonumber
&& + \frac{\mu^{\{k\}}}{2} \left\Vert \R_i^{-1/2} [ \mathcal{H}(\x) - \y_i ] - \z^{\{k\}}_i
- \frac{\bm{\lambda}_i^{\{k\}}}{\mu^{\{k\}} } \right\Vert_2^2.
\end{eqnarray}
%
%
%
To solve \eqref{eqn:L1_modified_3dvar} we carry out a regular $\Ltwo$-3D-Var minimization of the form \eqref{eqn:3DVar}
\begin{eqnarray}
\x^{\{k+1\}}_i &:=&
\arg\min_{\x}~ \frac{1}{2}
\left\Vert  \x - \x_i^{\rm b} \right\Vert_{\B_i^{-1}}^2 +
\frac{1}{2} \left\Vert \mathcal{H}(\x) - \y^{\{k\}}_i  \right\Vert_{ \R^{\{k\}\,-1}_i}^2,
\end{eqnarray}
but with  modified observations and a scaled covariance matrix:
\begin{equation}
\label{eqn:modified_data}
\y^{\{k\}}_i := \y_i + \R_i^{1/2} \, \left( \z^{\{k\}}_i + \frac{\bm{\lambda}_i^{\{k\}}}{\mu^{\{k\}} } \right),
\quad \R^{\{k\}}_i := \R_i/\mu^{\{k\}}, \\
\end{equation}
\end{subequations}
\item[Step 4:] Fix $\x_i^{\{k+1\}}$, $\bm{\lambda}_i^{\{k\}}$, and $\mu^{\{k\}}$, and solve
\begin{subequations}
\label{eqn:3dvar_L1_z_optimization}
\begin{equation}
\label{eqn:3dvar_L1_z_problem}
\z^{\{k+1\}}_i := \arg\min_\z \quad \| \z \|_{1} + \frac{\mu^{\{k\}}}{2}
\left\| \mathbf{d}^{\{k+1\}}_i - \z - \frac{\bm{\lambda}_i^{\{k\}}}{ \mu^{\{k\}} }\right\|_2^2,
\end{equation}
where
\begin{equation}
\label{eqn:increment}
\mathbf{d}^{\{k+1\}}_i :=  \R_i^{-1/2} \, \left[\mathcal{H}\left(\x^{\{k+1\}}_i\right) -\y_i\right].
\end{equation}
The above minimization subproblem can be solved by using the shrinkage procedure defined in Algorithm \ref{alg:L1Shrinkage}:
\begin{equation}
\label{eqn:shrinkage-L1}
\z^{\{k+1\}}_i := \textsc{LONEShrinkage}\left(\, \mu^{\{k\}};\, \mathbf{d}^{\{k+1\}}_i; \, \bm{\lambda}_i^{\{k\}}\, \right).
\end{equation}
\end{subequations}
\item[Step 5:] Update $\bm{\lambda}_i$:
\[
\bm{\lambda}_i^{\{k+1\}} := \bm{\lambda}_i^{\{k\}} - \mathbf{d}^{\{k+1\}}_i + \z^{\{k+1\}}_i.
\]
\item[Step 6:] Update $\mu$:
\[
\mu^{\{k+1\}} := \rho\, \mu^{\{k\}} , \qquad \rho >1.
\]
\end{enumerate}
\begin{algorithm}
\caption{$\Lone$\_Shrinkage}\label{alg:L1Shrinkage}
\begin{algorithmic}[1]
\Procedure{$\z$=LONEShrinkage($\mu; \mathbf{d};  \bm{\lambda}$)}{}
\State \textbf{Input:} $\lbrack \mu\in \mathbb{R};~ \mathbf{d}\in \mathbb{R}^m; ~ \bm{\lambda}\in \mathbb{R}^m\rbrack$
\State \textbf{Output:} $\lbrack \textbf{z} \in \mathbb{R}^m \rbrack$
\State{$\displaystyle
\z := \max \bigl\{ \,\left\| \mathbf{d} - \bm{\lambda}/ \mu  \right\|_2 - 1/\mu\, ,\, 0 \, \bigr\} \cdot
\frac{ \mathbf{d} - \bm{\lambda} / \mu }
{ \| \mathbf{d} - \bm{\lambda} / \mu \|_2 }.
$}
\EndProcedure
\end{algorithmic}
\end{algorithm}

%
%

\subsection{Huber-norm 3D-Var}
Using the $\Lone$ norm throughout spoils the smoothness properties near the mean. The smoothness property of $\Ltwo$ norm near the mean is highly desirable and can be retained by using Huber norm. The Huber norm treats the errors using $\Ltwo$ norm in the vicinity of the mean, and using $\Lone$ norm far from the mean \cite{lorenc1993bayesian}. From a statistical perspective, small scaled observation errors \eqref{eqn:define-z} are assumed to have independent standard normal distributions, while large scaled observation errors are assumed to have independent Laplace distributions \eqref{eqn:laplace-product-distribution}. The good properties of the traditional 3D-Var are preserved when the data is reasonable, but outliers are treated with the more robust $\Lone$ norm.

The 3D-Var problem in the Huber norm reads:
\begin{eqnarray}
\label{eqn:Huber-L1-ADMM}
&& \min\, \mathcal{J}(\x_i,\z_i) = \frac{1}{2}
\| \x_i - \x_i^{\rm b} \|_{\B_i^{-1}}^2 +
\frac{1}{2} \| \z_i\|_\textsc{hub} \\
\nonumber
&& \textnormal{subject to} \quad \z_i = \R_i^{-1/2}\, [ \mathcal{H}(\x_i) - \y_i ],
\end{eqnarray}
where
\begin{equation}
\label{eqn:huber-norm-definition}
\| \z_i \|_\textsc{hub} = \sum_{\ell=1}^m\, g_\textsc{hub}(\z_{i,\ell}), \qquad
g_\textsc{hub}(a) = \left\{
\begin{aligned}
a^2/2, & \quad |a| \le \tau \\
|a| - 1/2, & \quad |a| > \tau.
\end{aligned}
\right.
\end{equation}
Here $\z_{i,\ell}$ is the $\ell$-th entry of the vector $\z_i$. The threshold $\tau$ represents the number of standard deviations after which we switch from $\Ltwo$ to $\Lone$ norm.

\subsubsection{ADMM solution for Huber-3D-Var}

The Huber norm minimization problem \eqref{eqn:Huber-L1-ADMM} can be solved using the ADMM framework in a manner similar to the one described in Section \ref{sec:L13DVar}. The only difference is in Step 3: after fixing $\x_i^{\{k+1\}}$, $\bm{\lambda}_i^{\{k\}}$, and $\mu^{\{k\}}$ we do not solve the $\Lone$-norm problem \eqref{eqn:3dvar_L1_z_problem}. Rather, we solve the Huber-norm problem
\begin{subequations}
\label{eqn:3dvar_Huber_z_optimization}
\begin{equation}
\label{eqn:3dvar_Huber_z_problem}
\z^{\{k+1\}}_i := \arg\min_\z \quad \| \z \|_\textsc{hub} + \frac{\mu^{\{k\}}}{2}
\left\| \mathbf{d}^{\{k+1\}}_i - \z - \frac{\bm{\lambda}_i^{\{k\}}}{ \mu^{\{k\}} }\right\|_2^2,
\end{equation}
where $\mathbf{d}^{\{k+1\}}_i$ is defined in \eqref{eqn:increment}.
The closed-form solution of \eqref{eqn:3dvar_Huber_z_problem} is given by the Huber shrinkage procedure defined in Algorithm \ref{alg:HuberShrinkage}:
\[
\z^{\{k+1\}}_{i} := \textsc{HuberShrinkage}\left( \, \mu^{\{k\}};\, \mathbf{d}^{\{k+1\}}_i;\, \bm{\lambda}_i^{\{k\}}\,\right).
\]
\end{subequations}
\begin{algorithm}
\caption{Huber\_Shrinkage}\label{alg:HuberShrinkage}
\begin{algorithmic}[1]
\Procedure{$\z$=HuberShrinkage($\mu; \mathbf{d};  \bm{\lambda}$)}{}
\State \textbf{Input:} $\lbrack \mu\in \mathbb{R};~ \mathbf{d}\in \mathbb{R}^m; ~ \bm{\lambda}\in \mathbb{R}^m\rbrack$
\State \textbf{Output:} $\lbrack \textbf{z} \in \mathbb{R}^m \rbrack$
\For{$\ell=1,2,\dots,m$}
\If{$\lvert \mathbf{d}_\ell \rvert \ge \tau$}
\State{$\z_\ell = \displaystyle\max \left \{ \| \mathbf{d} - \bm{\lambda} / \mu \|_2 - 1/\mu, 0 \right \} \cdot
\frac{ \mathbf{d}_\ell - \bm{\lambda}_\ell / \mu }
{ \| \mathbf{d} - \bm{\lambda} / \mu \|_2 }.$}
\Else\State{$\z_\ell = \displaystyle\frac{\mu}{1+ \mu} ( \mathbf{d}_\ell - \bm{\lambda}_\ell / \mu )$}
\EndIf
\EndFor
\EndProcedure
\end{algorithmic}
\end{algorithm}

%
%

\subsubsection{Half-quadratic optimization solution for Huber-3D-Var}
\label{sec:3dvar-huber-halfq}

An alternative approach to carry out the minimization \eqref{eqn:Huber-L1-ADMM} is to reformulate it as a multiplicative half-quadratic optimization problem \cite{nikolova2005analysis}. We construct an augmented cost function which involves the auxiliary variable $\u_i \in \mathbb{R}^m$
\begin{subequations}
\label{eqn:Huber-A-function}
\begin{equation}
\mathcal{A}(\x_i,\u_i) =
\frac{1}{2}
\| \x_i - \x_i^{\rm b} \|_{\B_i^{-1}}^2 +
\frac{1}{2} \sum_{\ell=1}^{m} \bigl( \u_{i,\ell} \, \z_{i,\ell}^2 / 2 + \psi (\u_{i,\ell}) \bigr),
\end{equation}
where $\psi$ is a dual potential function determined by using the theory of convex conjugacy such that
\begin{equation}
\mathcal{J}(\x_i) = \min_{\u_i}\, \mathcal{A}(\x_i,\u_i).
\end{equation}
\end{subequations}
The minimizer of $\mathcal{A}$ is calculated using alternate minimization. Let the solution at the end of iteration $\{k-1\}$ read $(\x_i^{\{k\}}, \u_i^{\{k\}})$. At iteration $k$ we proceed as follows.
\begin{enumerate}
\item[Step 1:] First calculate $\u_i^{\{k+1\}}$ where $\x_i^{\{k\}}$ is fixed such, that
\[
\mathcal{A} \left(\x_i^{\{k\}}, \u_i^{\{k+1\}}\right) \le \mathcal{A}\left( \x_i^{\{k\}}, \u\right),
\quad \forall\, \u\,.
\]
This amounts to finding $\u_i^{\{k+1\}}$ according to
\begin{subequations}
\label{eqn:hqo-u}
\begin{equation}
\label{eqn:u-huber}
\u_i^{\{k+1\}} = \sigma \left( \R_i^{-1/2}\, \left[ \mathcal{H}(\x_i^{\{k\}}) - \y_i \right] \right)
\end{equation}
(where the function $\sigma$ is applied to each element of the argument vector).
For Huber regularization
\begin{equation}
\label{eqn:sigma-huber}
\sigma(a) = \left\{
\begin{array}{ll}
1,  &\quad |a| \le \tau; \\
\tau/|a|, &\quad |a| > \tau.
\end{array}
\right.
\end{equation}
\end{subequations}
\item[Step 2:] Next, calculate $\x_i^{\{k+1\}}$ where $\u_i^{(\{k+1\}}$ is fixed, such that
\[
\mathcal{A} \left(\x_i^{\{k+1\}}, \u_i^{\{k+1\}}\right) \le \mathcal{A} \left(\x, \u_i^{\{k+1\}}\right),
\quad \forall\, \x\,.
\]
This is achieved by solving the following minimization problem:
\begin{eqnarray*}
\x_i^{\{k+1\}}&=& \arg\min_{\x_i}~~
\frac{1}{2}
\| \x_i - \x_i^{\rm b} \|_{\B_i^{-1}}^2 +
\frac{1}{2} \sum_{\ell=1}^m \bigl( \u_{i,\ell}^{\{k+1\}}\, \z_{i,\ell}^2 / 2 \bigr) \\
&=& \arg\min_{\x_i}~~
\frac{1}{2}
\| \x_i - \x_i^{\rm b} \|_{\B_i^{-1}}^2 \\
&& +
\frac{1}{2} \left( \R_i^{-1/2}\, [ \mathcal{H}(\x_i) - \y_i ]\right )^T \, {\rm diag} \left( \u_i^{\{k+1\}}/2 \right)\,
\left( \R_i^{-1/2}\, [ \mathcal{H}(\x_i) - \y_i ] \right).
\end{eqnarray*}
This minimization problem is equivalent to solving a regular $\Ltwo$-3D-Var problem
\begin{equation}
\x_i^{\{k+1\}} = \arg\min_{\x_i}~~
\frac{1}{2} \| \x_i - \x_i^{\rm b} \|_{\B_i^{-1}}^2 + \| \mathcal{H}(\x_i) - \y_i \|_{\R_i^{\{k+1\}\, -1}}^2,
\end{equation}
with a modified observation error covariance
%
%
\begin{equation}
\label{eqn:huber-halfq-covariance}
\R_i^{\{k+1\}} := \R_i^{1/2}\cdot {\rm diag} \left( 2/\u_i^{\{k+1\}} \right) \cdot \R_i^{1/2}.
\end{equation}
For the case where $\R_i$ is diagonal equation \eqref{eqn:huber-halfq-covariance} divides each variance by $\u_i^{\{k+1\}}$. The smaller the entry $\u_{i,\ell}^{\{k+1\}}$ is the lower the weight given to observation $\ell$ in the 3D-Var data assimilation procedure becomes.

\end{enumerate}
%

\section{Robust 4D-Var}\label{sec:Robust-4DVar}
\subsection{$\Lone$-4D-Var}
\label{sec:L1-4DVar}

Given the background value of the initial state $\x_0^{\rm b}$, the covariance of the initial background errors $\B_0$, the observations $\y_i$ at $t_i$, and the corresponding observation error covariances $\R_i$ ($i=1,2,\cdots,N$), the $\Lone$-4D-Var problem looks for the MAP estimate $\x_0$ of the true initial conditions by solving the following constrained optimization problem:
\begin{equation}
\label{eqn:4dvar-L1-ADMM}
\begin{split}
\min_{\x_0}~\mathcal{J}(\x_0) & ~ := \frac{1}{2}
\| \x_0 - \x_0^{\rm b} \|_{\B_0^{-1}}^2 +
\frac{1}{2}
\sum_{i=1}^{N} \| \z_i  \|_1 \\
\textnormal{subject to:} & ~~\z_i = \R_i^{-1/2} [ \mathcal{H}(\x_i) - \y_i ] , \quad i = 1, 2, \cdots, N, \\
& ~~\x_{i} = \mathcal{M}_{i-1,i} (\x_{i-1}) , \quad i = 1, 2, \cdots, N.
\end{split}
\end{equation}
The augmented Lagrangian  associated with \eqref{eqn:4dvar-L1-ADMM} reads:
\begin{eqnarray}
\label{eqn:4dvar-Lagrangian-simple}
\qquad \mathcal{L} &=& \frac{1}{2} \left\Vert\x_0 - \x_0^{\rm b} \right\Vert^2_{ \mathbf{B}_0^{-1}}
 + \frac{1}{2} \sum_{\rm i=1}^{\rm N}\left\Vert \z_i \right\Vert_1 \\
 \nonumber
 && -\sum_{\rm i=1}^{\rm N} \theta_{i}^T\, \bigl(\x_i -\mathcal{M}_{i-1, i}(\x_{i-1})\bigr)
 + \frac{\nu}{2}\, \sum_{i=1}^{\rm N} \left\Vert \x_i - \mathcal{M}_{i-1, i}(\x_{i-1}) \right\Vert^2_{\mathbf{P}_{i}^{-1}}\\
\nonumber
 && -\sum_{\rm i=1}^{\rm N} \bm{\lambda}_{i}^T\, \left(\R_i^{-1/2} [ \mathcal{H}(\x_i) - \y_i ] - \z_i \right)
 + \frac{\mu}{2}\, \sum_{i=1}^{\rm N} \left\Vert \R_i^{-1/2} [ \mathcal{H}(\x_i) - \y_i ] - \z_i \right\Vert^2_2\,,
\end{eqnarray}
where $\mathbf{P}_{i}$'s are model error scaling matrices.

Using the identity \eqref{eqn:z-identity} the Lagrangian \eqref{eqn:4dvar-Lagrangian-simple} becomes
\begin{eqnarray}
\label{eqn:4dvar-L1-Lagrangian}
 \mathcal{L} &=& \frac{1}{2} \left\Vert\x_0 - \x_0^{\rm b} \right\Vert^2_{ \mathbf{B}_0^{-1}}
 + \frac{1}{2} \sum_{\rm i=1}^{\rm N}\left\Vert \z_i \right\Vert_1 \\
 \nonumber
 && -\sum_{\rm i=1}^{\rm N} \theta_{i}^T\, \bigl(\x_i -\mathcal{M}_{i-1, i}(\x_{i-1})\bigr)
   + \frac{\nu}{2}\, \sum_{i=1}^{\rm N} \left\Vert \x_i - \mathcal{M}_{i-1, i}(\x_{i-1}) \right\Vert^2_{\mathbf{P}_{i}^{-1}}\\
\nonumber
 && + \frac{\mu}{2}\,\sum_{\rm i=1}^{\rm N}\,\left\Vert \R_i^{-1/2} [ \mathcal{H}(\x_i) - \y_i ] - \z_i - \bm{\lambda}_i/ \mu \right\Vert_2^2
 - \frac{1}{2\mu}\,\sum_{\rm i=1}^{\rm N}\, \left\Vert  \bm{\lambda}_i \right\Vert_2^2
\end{eqnarray}

The problem \eqref{eqn:4dvar-L1-ADMM} can be solved by alternating minimization of the Lagrangian \eqref{eqn:4dvar-L1-Lagrangian}, as follows:
\begin{enumerate}
\item[Step 1:] Initialize $\x_0^{\{0\}} = \x_0^{\rm b}\,$.
\item[Step 2:] Run the forward model starting from $\x^{\{0\}}_0$ to obtain $\x^{\{0\}}_i$, $i=1,\dots,N$.
\item[Step 3:] Set $\z_i^{\{0\}} = \R_i^{-1/2}\, [ \mathcal{H}(\x_i^{\{0\}}) - \y_i ]\,$ , $\bm{\lambda}_i^{\{0\}} = \bm{0}\,$, for $i=1,\dots,N$,  $\mu^{\{0\}} = 1$.
\item[Step 4:]
Fix $\z_i^{\{k\}}$ and $\bm{\lambda}_i^{\{k\}}$, and minimize the Lagrangian to solve for $\x_0^{\{k+1\}}$:
\begin{eqnarray*}
\x_0^{\{k+1\}} &:=& \arg\min~ \mathcal{L},\\
\mathcal{L} &=& \frac{1}{2} \left\Vert\x_0 - \x_0^{\rm b} \right\Vert^2_{ \mathbf{B}_0^{-1}}
-\sum_{\rm i=1}^{\rm N} \theta_{i}^T\, \bigl(\x_i -\mathcal{M}_{i-1, i}(\x_{i-1})\bigr) \\
 &&\nonumber  + \frac{\nu}{2}\, \sum_{i=1}^{\rm N} \left\Vert \x_i - \mathcal{M}_{i-1, i}(\x_{i-1}) \right\Vert^2_{\mathbf{P}_{i}^{-1}}\\
\nonumber
 && + \frac{\mu^{\{k\}}}{2}\,\sum_{\rm i=1}^{\rm N}\,\left\Vert \R_i^{-1/2} [ \mathcal{H}(\x_i) - \y_i ] - \z_i^{\{k\}} - \bm{\lambda}_i^{\{k\}}/ \mu^{\{k\}} \right\Vert_2^2
\end{eqnarray*}
This is equivalent to solving the traditional $\Ltwo$-4D-Var problem
\begin{eqnarray*}
\x_0^{\{k+1\}} := \arg\min ~\mathcal{J}(\x_0) &=& \frac{1}{2}
\| \x_0 - \x_0^{\rm b} \|_{\B_0^{-1}}^2 +
\frac{1}{2}
\sum_{i=1}^{N} \|  \mathcal{H}(\x_i) - \widetilde{\y}_i \|_{\widetilde{\R}_i^{-1}}^2 \\
\textnormal{subject to} \quad
\x_{i+1} &=& \mathcal{M}_{i,i+1} (\x_i) , \quad i = 0, 1, 2, \cdots, N -1,
\end{eqnarray*}
with the modified observation vectors and observation error covariances \eqref{eqn:modified_data}.
\item[Step 5:]
Run the forward model starting from $\x_0^{\{k+1\}}$ to obtain $\x_i^{\{k+1\}}$, $i=1,\dots,N$.
\item[Step 6:]
Fix $\x_i^{\{k+1\}}$ and $\bm{\lambda}_i^{\{k\}}$, and find $\z_i^{\{k+1\}}$, $i=1,\dots,N$ as
\begin{eqnarray}
\nonumber
&&\qquad \mathbf{z}^{\{k+1\}} := \arg\min~ \mathcal{L}, \\
\label{eqn:4dvar-z}
&&\qquad \mathcal{L} =  \frac{1}{2} \sum_{\rm i=1}^{\rm N}\left\Vert \z_i \right\Vert_1+ \frac{\mu}{2}\,\sum_{\rm i=1}^{\rm N}\,\left\Vert \R_i^{-1/2} \mathcal{H}(\x_i^{\{k+1\}}) - \y_i ] - \z_i - \bm{\lambda}_i^{\{k\}}/ \mu^{\{k\}} \right\Vert_2^2 \\
\nonumber
&& \qquad\qquad = \frac{1}{2} \left\Vert \z \right\Vert_1
+ \frac{\mu}{2}\,\left\Vert \bm{d}^{\{k+1\}}  - \bm{\lambda}^{\{k\}}/ \mu^{\{k\}} - \z \right\Vert_2^2,
\end{eqnarray}
where we denote
\[
\z = \begin{bmatrix} \z_1 \\ \vdots \\ \z_N \end{bmatrix}, \quad
\bm{\bm{\lambda}}^{\{k\}} = \begin{bmatrix} \bm{\lambda}_1^{\{k\}} \\ \vdots \\ \bm{\lambda}_N^{\{k\}} \end{bmatrix}, \quad
\bm{d}^{\{k+1\}} = \begin{bmatrix} \R_1^{-1/2} [ \mathcal{H}(\x_1^{\{k+1\}}) - \y_1 ]  \\ \vdots \\ \R_N^{-1/2} [ \mathcal{H}(\x_N^{\{k+1\}}) - \y_N ]  \end{bmatrix}.
\]
Problem \eqref{eqn:4dvar-z}  is solved using  Algorithm \ref{alg:L1Shrinkage}:
\begin{equation}
\label{eqn:shrinkage-L1-4dvar}
\z^{\{k+1\}} := \textsc{LoneShrinkage}\left(\, \mu^{\{k\}};\, \bm{d}^{\{k+1\}}; \, \bm{\lambda}^{\{k\}}\, \right).
\end{equation}

\item[Step 7:]
Update $\bm{\bm{\lambda}}$:
\[
\bm{\lambda}^{\{k+1\}}_i := \bm{\lambda}^{\{k\}}_i - \R_i^{-1/2} \, \left[ \mathcal{H}\left(\x^{\{k+1\}}_i\right) - \y_i \right] + \z^{\{k+1\}}_i, \quad i = 1,\dots,N.
\]
\end{enumerate}

\subsection{Huber 4D-Var}
\label{sec:Huber-4DVar}

The Huber 4D-Var problem reads
\begin{equation}
\label{eqn:4dvar-Huber-ADMM}
\begin{split}
\min_{\x_0} ~ \mathcal{J}(\x_0,\z) &~ := \frac{1}{2}
\| \x_0 - \x_0^{\rm b} \|_{\B_0^{-1}}^2 +
\frac{1}{2}
\sum_{i=1}^{N} \| \z_i  \|_{\textsc{hub}} \\
\textnormal{subject to:} ~~ \z_i &~= \R_i^{-1/2} [ \mathcal{H}(\x_i) - \y_i ] , \quad i = 1, 2, \cdots, N, \\
\x_{i} &~= \mathcal{M}_{i-1,i} (\x_{i-1}) , \quad i = 1, 2, \cdots, N.
\end{split}
\end{equation}
%

\subsubsection{ADMM solution of Huber-4D-Var}
The Huber norm minimization problem \eqref{eqn:4dvar-Huber-ADMM} can be solved using the ADMM framework in a manner similar to the one described in Section \ref{sec:L1-4DVar}. The only difference is in step 3, where we find $\z$ by solving the Huber-norm problem
\begin{eqnarray}
\label{eqn:4dvar-Huber-z-alt}
\z^{\{k+1\}} := \arg\min_\z~ \mathcal{L} &=&  \frac{1}{2} \left\Vert \z \right\Vert_\textsc{hub}
+ \frac{\mu^{\{k\}}}{2}\,\left\Vert \bm{d}^{\{k+1\}}  - \bm{\lambda}^{\{k\}}/ \mu^{\{k\}} - \z \right\Vert_2^2\,.
\end{eqnarray}
The closed form solution of \eqref{eqn:4dvar-Huber-z-alt} is obtained by the Huber shrinkage procedure defined in Algorithm \ref{alg:HuberShrinkage}:
\[
\z^{\{k+1\}} := \textsc{HuberShrinkage}\left( \, \mu^{\{k\}};\, \mathbf{d}^{\{k+1\}};\, \bm{\lambda}^{\{k\}}\,\right).
\]

\subsubsection{Half-quadratic optimization solution for Huber-4D-Var}
An alternative procedure to solve the 4D-Var problem described in \eqref{eqn:4dvar-Huber-ADMM} is using the steps described in Section \ref{sec:3dvar-huber-halfq}.  Use \eqref{eqn:hqo-u} to update the scaling factors $\u_i^{\{k+1\}}$ for each $i=1,\dots,N$. Use   \eqref{eqn:huber-halfq-covariance} to compute the scaled observation error covariance matrices  $\R_i^{\{k+1\}}$ for each $i=1,\dots,N$. Solve a regular $\Ltwo$-4D-Var problem \eqref{eqn:L2-4dvar} with the modified observation error covariance matrices to obtain the analysis $\x_0^{\{k+1\}}$. Propagate the initial solution in time to obtain the forward solution $\x_i^{\{k+1\}}$, $i=1,\dots,N$. Iterate by computing again the new scaling factors $\u_i^{\{k+2\}}$ and so on.

\section{Robust data assimilation using the ensemble Kalman filter} \label{sec:L1-EnKF}

Robust EnKF algorithms are obtained by reformulating the equivalent optimization problems \eqref{eqn:EnKF-L2} or \eqref{eqn:enkf-traditional-opt-L2} using  $\| \mathcal{H}\left(\overline{\x}^{\rm b}_i + \mathbf{X}^{\rm b}_i \, w\right) - \y_i \|_*$, where $\| \cdot \|_*$ stands for either $\Lone$ or Huber norms. The solution of the modified optimization problems relies on repeated applications of the standard EnKF procedure.
%
%
%
%
For traditional EnKF \eqref{eqn:enkf-traditional-opt-L2} a different modified optimization problem has to be solved for each ensemble member. This is computationally expensive. For this reason we discuss below only the robust ensemble square root filter.

\subsection{$\Lone$-EnKF}

The $\Lone$-EnKF problem reads:
\begin{equation}
\label{eqn:EnKF-L1-ADMM}
\begin{split}
&\min\, \mathcal{J}(w_i,\z_i) = (\nens-1)\, \| w \|_2^2 + \| \z_i \|_{1} \\
& \textnormal{subject to} \quad
\z_i = \R_i^{-1/2}\, \left[ \mathcal{H}\left(\overline{\x}^{\rm b}_i + \mathbf{X}^{\rm b}_i \, w_i\right) - \y_i  \right].
\end{split}
\end{equation}

The augmented Lagrangian for problem \eqref{eqn:EnKF-L1-ADMM} is given by
\begin{eqnarray}
\label{eqn:3dvar-augmented-Lagrangian-2}
\mathcal{L} &=& (\nens-1)\, \| w_i \|_2^2 +  \| \z_i \|_{1} \\
\nonumber
&& - \bm{\lambda}^T\, \left(  \R_i^{-1/2}\, \left[ \mathcal{H}\left(\overline{\x}^{\rm b}_i + \mathbf{X}^{\rm b}_i \, w_i\right) - \y_i  \right] - \z_i  \right) \\
\nonumber
&& + \frac{\mu}{2}\, \| \R_i^{-1/2}\, \left[ \mathcal{H}\left(\overline{\x}^{\rm b}_i + \mathbf{X}^{\rm b}_i \, w_i\right) - \y_i  \right] - \z_i \|_2^2 \\
\nonumber
&=& (\nens-1)\, \| w_i \|_2^2 +  \| \z_i \|_{1}  - \frac{1}{2\mu}\, \left\Vert  \bm{\lambda} \right\Vert_2^2 \\
\nonumber
&& + \frac{\mu}{2}\,\left\Vert \R_i^{-1/2}\, \left[ \mathcal{H}\left(\overline{\x}^{\rm b}_i + \mathbf{X}^{\rm b}_i \, w_i\right) - \y_i  \right] - \z_i - \bm{\lambda}/ \mu \right\Vert_2^2.
\end{eqnarray}

The problem \eqref{eqn:EnKF-L1-ADMM} is solved by ADMM as follows:
\begin{enumerate}
\item[Step 1:]
Fix $\mu^\pk$, $\z_i^\pk$ and $\boldsymbol{\lambda}_i^\pk$, and solve
\begin{equation}
\label{eqn:EnKF-L1-x}
\begin{split}
&~ w_i^\pkp := \arg\min_{w}~\mathcal{J}(w) =  (\nens-1)\, \| w \|_2^2 \\
&\qquad + \frac{\mu^\pk}{2}\,\left\Vert \R_i^{-1/2}\, \left[ \mathcal{H}\left(\overline{\x}^{\rm b}_i + \mathbf{X}^{\rm b}_i \, w\right) - \y_i  \right] - \z_i^\pk - \frac{\bm{\lambda}_i^\pk}{ \mu^\pk } \right\Vert_2^2.
\end{split}
\end{equation}
For a general nonlinear observation operator \eqref{eqn:EnKF-L1-x} can be solved by a regular nonlinear minimization in the ensemble space.
It is useful to compare \eqref{eqn:EnKF-L1-x} with the $\Ltwo$-EnKF equation \eqref{eqn:EnKF-L2}.
The weights update \eqref{eqn:EnKF-L1-x} can be obtained by applying the EnKF methodology to a data assimilation problem with the modified observations $\y^{\{k\}}_i$ and the scaled observation error covariance matrix $\R^{\{k\}}_i$ given by \eqref{eqn:modified_data}.
\item[Step 2:]
Fix $\mu^\pk$, $w_i^\pkp$ and $\bm{\lambda}_i^\pk$, and solve for $\z_i$
\[
\begin{split}
&~\z_i^\pkp := \arg\min_\z \, \mathcal{J}(\z), \\
&~\mathcal{J}(\z) = \| \z \|_{1} + \frac{\mu^\pk}{2} \left\| \R_i^{-1/2} \, \left[ \mathcal{H}\left(\overline{\x}^{\rm b}_i + \mathbf{X}^{\rm b}_i \, w_i^\pkp\right)  -\y_i\right] - \z - \frac{\bm{\lambda}^\pk }{ \mu^\pk } \right\|_2^2.
\end{split}
\]
The above minimization subproblem is solved by using the shrinkage procedure defined in Algorithm \ref{alg:L1Shrinkage}:
\begin{eqnarray}
\label{eqn:shrinkage-L1-EnKF}
\begin{split}
\mathbf{d}^{\{k+1\}}_i :=&~ \R_i^{-1/2}\, \left[ \mathcal{H}\left(\overline{\x}^{\rm b}_i + \mathbf{X}^{\rm b}_i \, w_i^\pkp\right)  -\y_i \right]\,  \\
\z^{\{k+1\}} :=&~ \textsc{LoneShrinkage}\left(\, \mu^{\{k\}};\, \mathbf{d}_i^{\{k+1\}}; \, \bm{\lambda}^{\{k\}}\, \right).
\end{split}
\end{eqnarray}
\item[Step 3:]
Update $\bm{\lambda}_i$:
\[
\bm{\lambda}_i^\pkp := \bm{\lambda}_i^\pk - \R_i^{-1/2} \, \left[ \mathcal{H}\left(\overline{\x}^{\rm b}_i + \mathbf{X}^{\rm b}_i \, w_i^\pkp\right) - \y_i \right] + \z_i^\pkp.
\]
\end{enumerate}

Assume that after $M$ iterations we achieve satisfactory convergence and stop. The mean analysis provided by the robust EnSRF is given by the optimal weights:
\begin{equation}
\label{eqn:robust-ensrf-mean}
\overline{w}_i^\textnormal{a} := w_i^{\{ M \}}.
\end{equation}
The weights of individual analysis ensemble are obtained by applying \eqref{eqn:letkf-analysis-ensemble} with the modified observation covariance of the last iteration:
\begin{equation}
\label{eqn:robust-ensrf-ensemble}
\begin{split}
&\left( (\nens-1) \mathbf{I} +  \mu^{\{ M \}}\,\left(\mathbf{Y}^{\rm b}_i\right)^T\, \R_i^{-1}\, \mathbf{Y}^{\rm b}_i \right)^{-1}
= (\nens-1)^{-1}\, \mathbf{W}_i^{\{ M \}}\, \mathbf{W}_i^{\{ M \}}\,^T, \\
&w_i^{\textnormal{a}\langle\ell\rangle} := \overline{w}^\textnormal{a}_i + \mathbf{W}_i^{\{ M \}}(:,\ell),
\quad
\x_i^{\textnormal{a}\langle\ell\rangle} = \overline{\x}_i^\textnormal{b} + \mathbf{X}_i^\textnormal{b} \, w_i^{\textnormal{a}\langle\ell\rangle}, \quad \ell = 1,\dots,\nens.
\end{split}
\end{equation}

\subsection{Huber EnKF}

The Huber EnKF problem reads:
\begin{equation}
\label{eqn:EnKF-Huber-HalfQuadratic}
\begin{split}
& \min\, J(w,\z_i) = (\nens-1)\, \| w \|_2^2 +  \| \z_i\|_\textsc{hub} \\
& \textnormal{subject to} \quad
\z_i = \R_i^{-1/2}\, \left[ \mathcal{H}\left(\overline{\x}^{\rm b}_i + \mathbf{X}^{\rm b}_i \, w\right) - \y_i  \right].
\end{split}
\end{equation}
The problem \eqref{eqn:EnKF-Huber-HalfQuadratic} is solved iteratively by half-quadratic minimization. The minimizer of $\mathcal{A}$ is calculated using alternate minimization. Let the solution at iteration $\{k\}$ be $w_i^{\{k\}}, \u_i^{\{k\}}$. We proceed as follows.
\begin{enumerate}
\item[Step 1:] First calculate
\[
\u_i^\pkp = \sigma \left( \R_i^{-1/2}\, \left[ \mathcal{H}\left(\overline{\x}^{\rm b}_i + \mathbf{X}^{\rm b}_i \, w_i^\pk\right) - \y_i \right] \right)
\]
where the function $\sigma$ is given by \eqref {eqn:sigma-huber}.

\item[Step 2:] Next, update the state $\x_i^\pkp = \overline{\x}_i^\textnormal{b} + \mathbf{X}_i^\textnormal{b}\, w_i^\pkp$ as follows:
\begin{eqnarray*}
w_i^\pkp &=& \arg\min_{w}~~
(\nens-1)\, \| w \|_2^2
 +
 \sum_\ell \left[ (\u_i^\pkp)_\ell\, z_\ell^2 / 2 \right] \\
&=& \arg\min_{w}~~
 (\nens-1)\, \| w \|_2^2 + \| \mathcal{H}\left(\overline{\x}^{\rm b}_i + \mathbf{X}^{\rm b}_i \, w\right) - \y_i \|_{\R_i^\pkp\,^{-1}}^2 \\
\textnormal{where}&& \R_i^\pkp\,^{-1} = \R_i^{-1/2}\, {\rm diag} \left( \u^{\{k\}}/2 \right) \,\R_i^{-1/2}.
\end{eqnarray*}
This second step requires the application of EnSRF \eqref{eqn:EnKF-L2} with modified observation covariance matrix. The analysis weights provide the new iterate $w_i^\pkp \equiv w_i^\textnormal{a}$.
\end{enumerate}

Assume that after $M$ iterations we achieve satisfactory convergence and stop. The mean analysis provided by the robust EnSRF is given by the optimal weights \eqref{eqn:robust-ensrf-mean}.
The weights of individual analysis ensemble members are obtained from \eqref{eqn:robust-ensrf-ensemble} with the modified observation covariance of the last iteration:
%
%
\[
\begin{split}
&\left( (\nens-1) \mathbf{I} +  \left(\mathbf{Y}^{\rm b}_i\right)^T\, \R_i^{-1/2}\, {\rm diag} \left( \frac{\u^{\{M\}}}{2} \right) \,\R_i^{-1/2}\, \mathbf{Y}^{\rm b}_i \right)^{-1} \\
& \qquad = (\nens-1)^{-1}\, \mathbf{W}_i^{\{ M \}}\, \mathbf{W}_i^{\{ M \}}\,^T.
\end{split}
\]

\section{Numerical experiments} \label{sec:numexp}

We compare the performance of variational and ensemble data assimilation techniques using $\Lone$, $\Ltwo$, and Huber norms using the Lorenz-96 model \cite{lorenz1996} and a shallow water model in spherical coordinates ($\sim$8,000 variables) \cite{Amik:2007}.  The error statistic used to compare analyses against the reference solution is the root mean squared error (RMSE) metric:
                  \begin{equation}
                  \label{eqn:RMSE_Formula}
                           \mathbf{RMSE}_i =
                           \nvar^{-1/2}\, \left\Vert \x_i - \x_i^{\rm true} \right\Vert_2
                  \end{equation}
                  where $\x^{\rm true}$ is the reference state of the system. The RMSE is calculated for each time within the assimilation window. We describe next the two testing models.

\subsection{Lorenz-96 model}
The Lorenz-96 model \cite{lorenz1996} is given by the equations
      \begin{equation}
      \label{eqn:Lorenz96}
               \frac{\mathrm{d}\x_{k}}{\mathrm{d}t} = \x_{\rm k-1} \left( \x_{\rm k+1} - \x_{\rm k-2} \right) - \x_{k} + F \,,
               \quad k=1,\dots,40,
      \end{equation}
with periodic boundary conditions and the forcing term $F=8$ \cite{lorenz1996}. A vector of equidistant components ranging from $-2$ to $2$ is integrated for one time unit, and the result is taken as the reference initial state for the experiments. The background errors at initial time are Gaussian, with a diagonal background error covariance matrix, and each component's standard deviation equal to $8\%$ of the time-averaged magnitude of the reference solution. The same covariance matrix is used for all times in the 3D-Var calculations. Synthetic observations  are generated by perturbing the reference trajectory with normal noise with mean zero, diagonal observation error covariance, and standard deviation of each component equal to $5\%$ of the time-averaged magnitude of the reference solution.

\subsection{Shallow water model on the sphere}
The shallow water equations have been used extensively as a simple model of the atmosphere since they contain the essential wave propagation mechanisms found in general circulation models  \cite{Amik:2007}. The shallow water equations in spherical coordinates are:
\begin{subequations}
\label{eqn:swe}
\begin{eqnarray}
 \frac{\partial u}{\partial t} + \frac{1}{a\cos \theta} \left( u \frac{\partial u}{\partial \bm{\lambda}} + v \cos \theta \frac{\partial u}{\partial \theta} \right) - \left(f + \frac{u \tan \theta}{a} \right) v + \frac{g}{a \cos \theta} \frac{\partial h} {\partial \bm{\lambda}} = 0, \\
 \frac{\partial v}{\partial t} + \frac{1}{a\cos \theta} \left( u \frac{\partial v}{\partial \bm{\lambda}} + v \cos \theta \frac{\partial v}{\partial \theta} \right) + \left(f + \frac{u \tan \theta}{a} \right) u + \frac{g}{a} \frac{\partial h} {\partial \theta} = 0, \\
 \frac{\partial h}{\partial t} + \frac{1}{a \cos \theta} \left(\frac{\partial\left(hu\right)}{\partial \bm{\lambda}} + \frac{\partial{\left(hv \cos \theta \right)}}{\partial \theta} \right) = 0.
\end{eqnarray}
\end{subequations}
Here $f$ is the Coriolis parameter given by $f = 2 \Omega \sin \theta$, where  $\Omega$ is the angular speed of the rotation of the Earth, $h$ is the height of the homogeneous atmosphere, $u$ and $v$ are the zonal and meridional wind components, respectively, $\theta$ and $\bm{\lambda}$ are the latitudinal and longitudinal directions, respectively, $a$ is the radius of the earth and $g$ is the gravitational constant. The space discretization is performed using the unstaggered Turkel-Zwas scheme \cite{Navon19911311, navon1987application}. The discretized spherical grid has nlon = 36 nodes in longitudinal direction and nlat = 72 nodes in the latitudinal direction. The semi-discretization in space leads to a discrete model of the form \eqref{eqn:genmodel}.
 In \eqref{eqn:genmodel} the zonal wind, meridional wind and the height variables are combined into the vector $\x \in \mathbb{R}^n$ with $n=3\times{\rm nlat}\times{\rm nlon}$. We perform the time integration using an adaptive time-stepping algorithm. For a tolerance of $\displaystyle 10^{-8}$ the average time step size of the time-integrator is 180 seconds. A reference initial condition is used to generate a reference trajectory.

Synthetic observation errors at various times $t_{i}$ are normally distributed  with mean zero and a diagonal observation error covariance matrix with entries equal to $(\mathbf{R}_{i})_{k,k}=1\,m^2 s^{-2}$ for $u$ and $v$ components, and
 $(\mathbf{R}_{i})_{k,k}=10^6\,m^2$ for $h$ components. These values correspond (approximately) to a standard deviation of $5\%$ of the time-averaged values for $u$ and $v$ components, and $2\%$ of the time-averaged values for $h$ component. A flow dependent background error covariance matrix is constructed as described in \cite{Attia:2014, Attia:2015}. The standard deviation of the background errors for the height component is 2\% of the average magnitude of the height component in the reference initial condition. The standard deviation of the background errors for the wind components is 15\% of the average magnitude of the wind component in the reference initial condition.
\subsection{3D-Var experiments with the Lorenz-96 model}

The assimilation window is 2 units long. All components of the state space are observed. Experiments are performed for two observation frequencies, namely, one observation for every 0.1 units and 0.01 units, respectively. The tests are carried out with both good and erroneous data. The good data contains no outliers. The erroneous data contains one outlier (measurement of one state) which lies $\sim 100$ standard deviations away from the mean. This outlier occurs every 0.2 time units. 

\begin{figure}[ht]
  \centering
  \subfigure[Observations with only small random errors \label{fig:3DVarRMSE_Lorenz_Good_0_01_1}]
  {\includegraphics[width=\myfigurewidth]{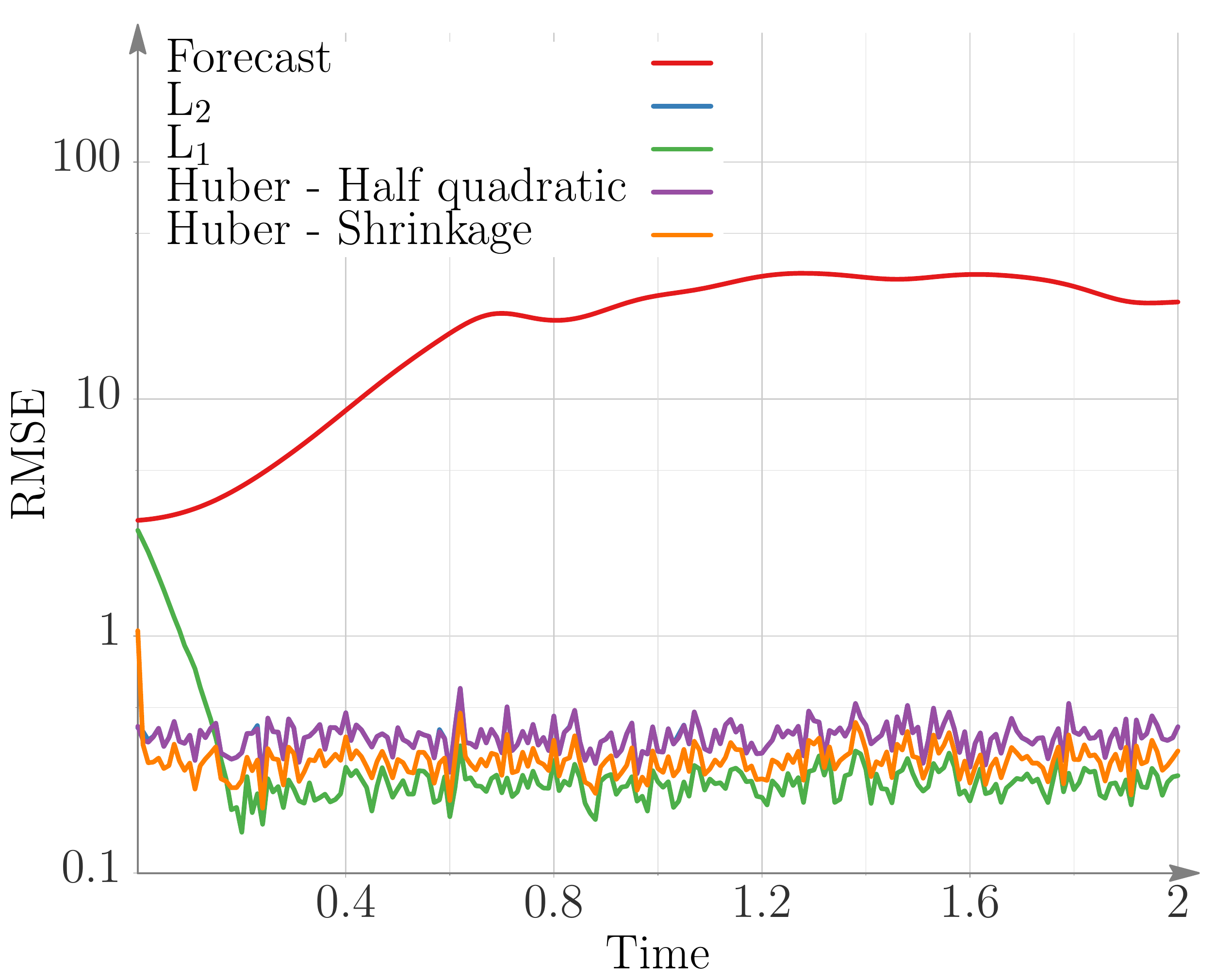}}
  \subfigure[Observations with outliers \label{fig:3DVarRMSE_Lorenz_Bad_0_01_1}]
  {\includegraphics[width=\myfigurewidth]{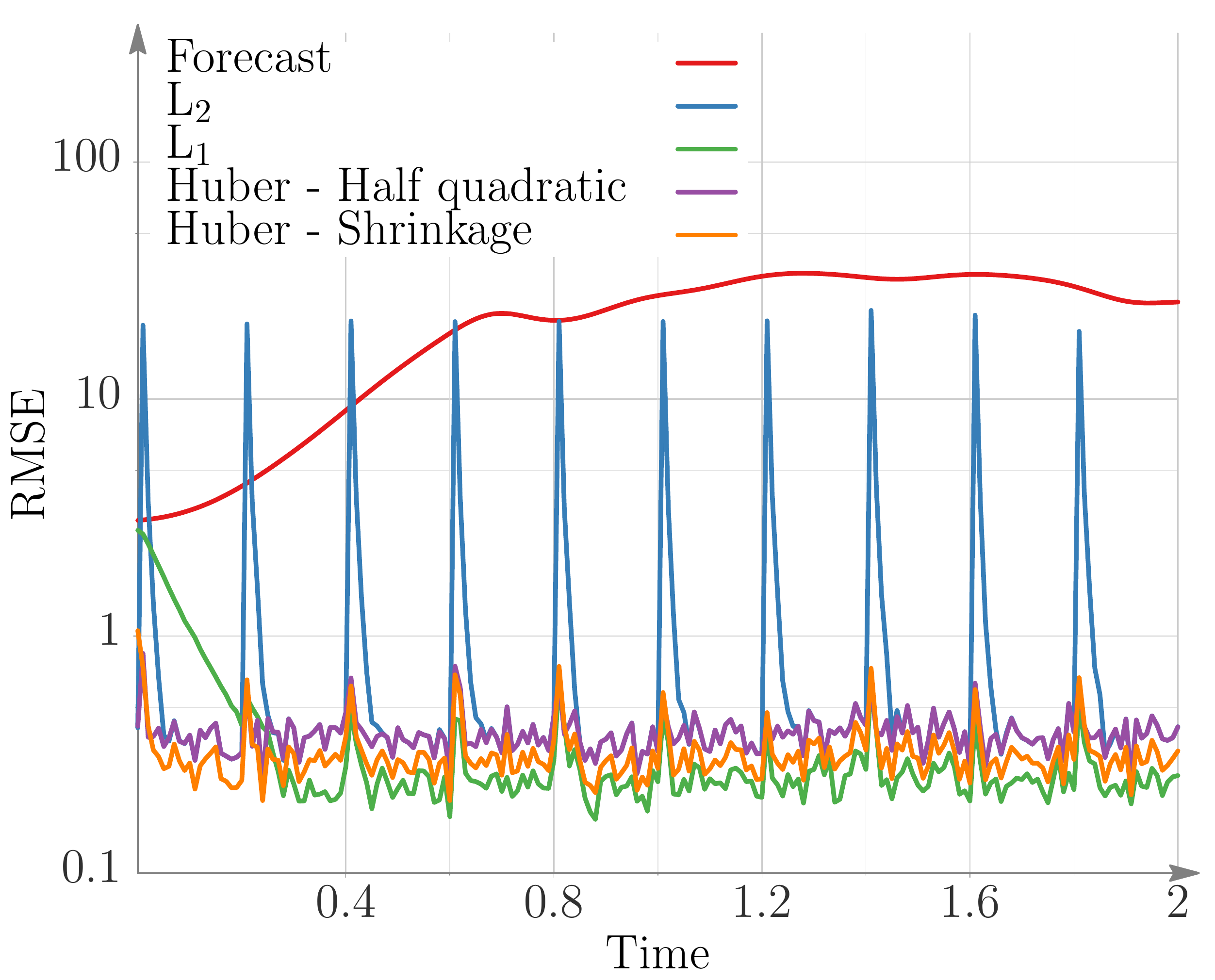}}
  \caption{3D-Var results for the Lorenz-96 model \eqref{eqn:Lorenz96}. The frequency of observations is 0.01 time units. Erroneous observations occur every 0.2 time units. The Huber norm uses $\tau = 1$.}
  \label{fig:Lorenz3DVar_0_01_1}
\end{figure}

\begin{figure}[ht]
  \centering
  \subfigure[Observations with only small random errors \label{fig:3DVarRMSE_Lorenz_Good_0_1_1}]
  {\includegraphics[width=\myfigurewidth]{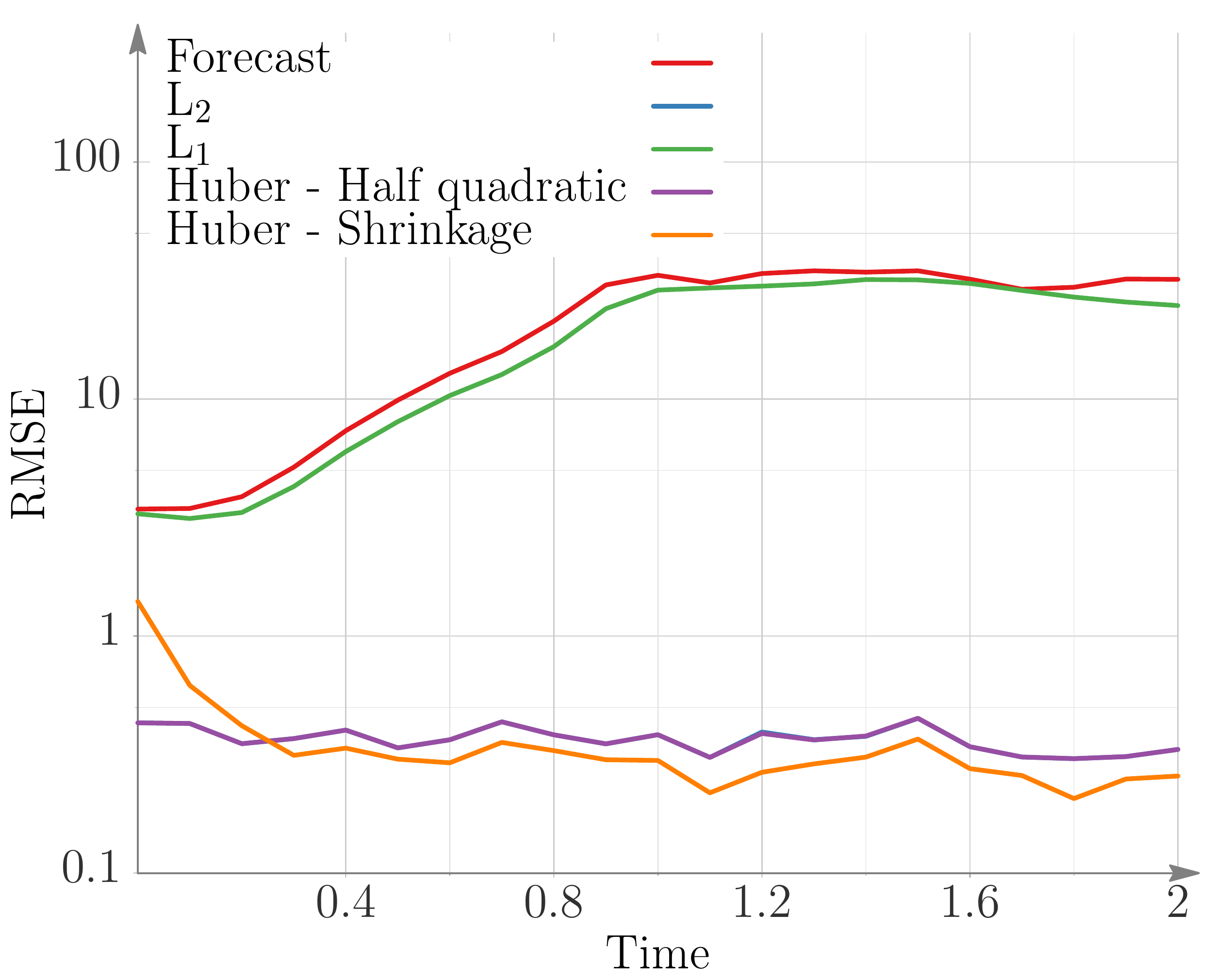}}
  \subfigure[Observations with outliers \label{fig:3DVarRMSE_Lorenz_Bad_0_1_1}]
  {\includegraphics[width=\myfigurewidth]{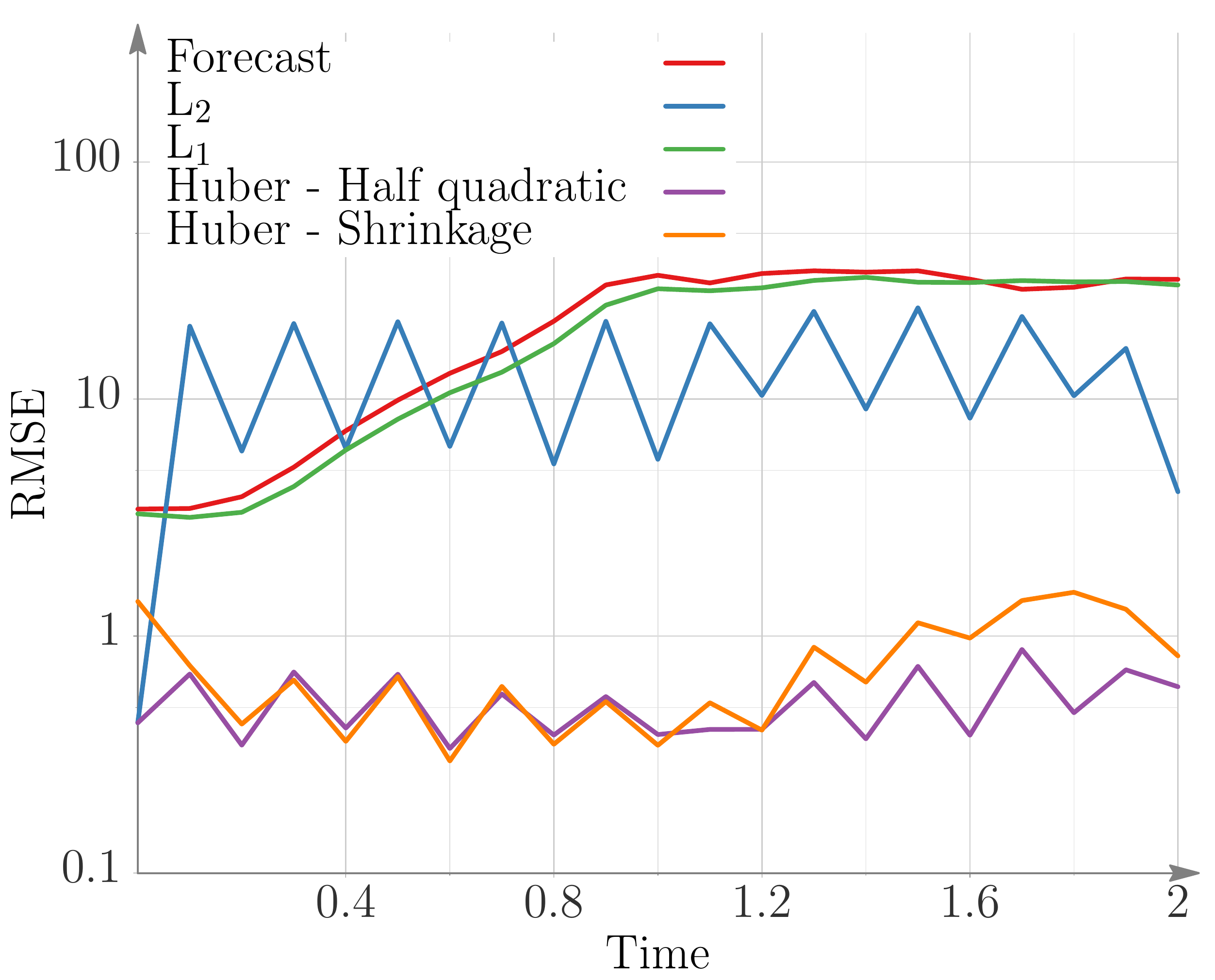}}
  \caption{3D-Var results for the Lorenz-96 model \eqref{eqn:Lorenz96}. The frequency of observations is 0.1 time units. Erroneous observations occur every 0.2 time units. The Huber norm uses $\tau = 1$.}
  \label{fig:Lorenz3DVar_0_1_1}
\end{figure}

\begin{figure}[ht]
  \centering
  \subfigure[Observations with only small random errors. Observations taken every 0.01 time units. \label{fig:3DVarRMSE_Lorenz_Good_0_01_3}]
  {\includegraphics[width=\myfigurewidth]{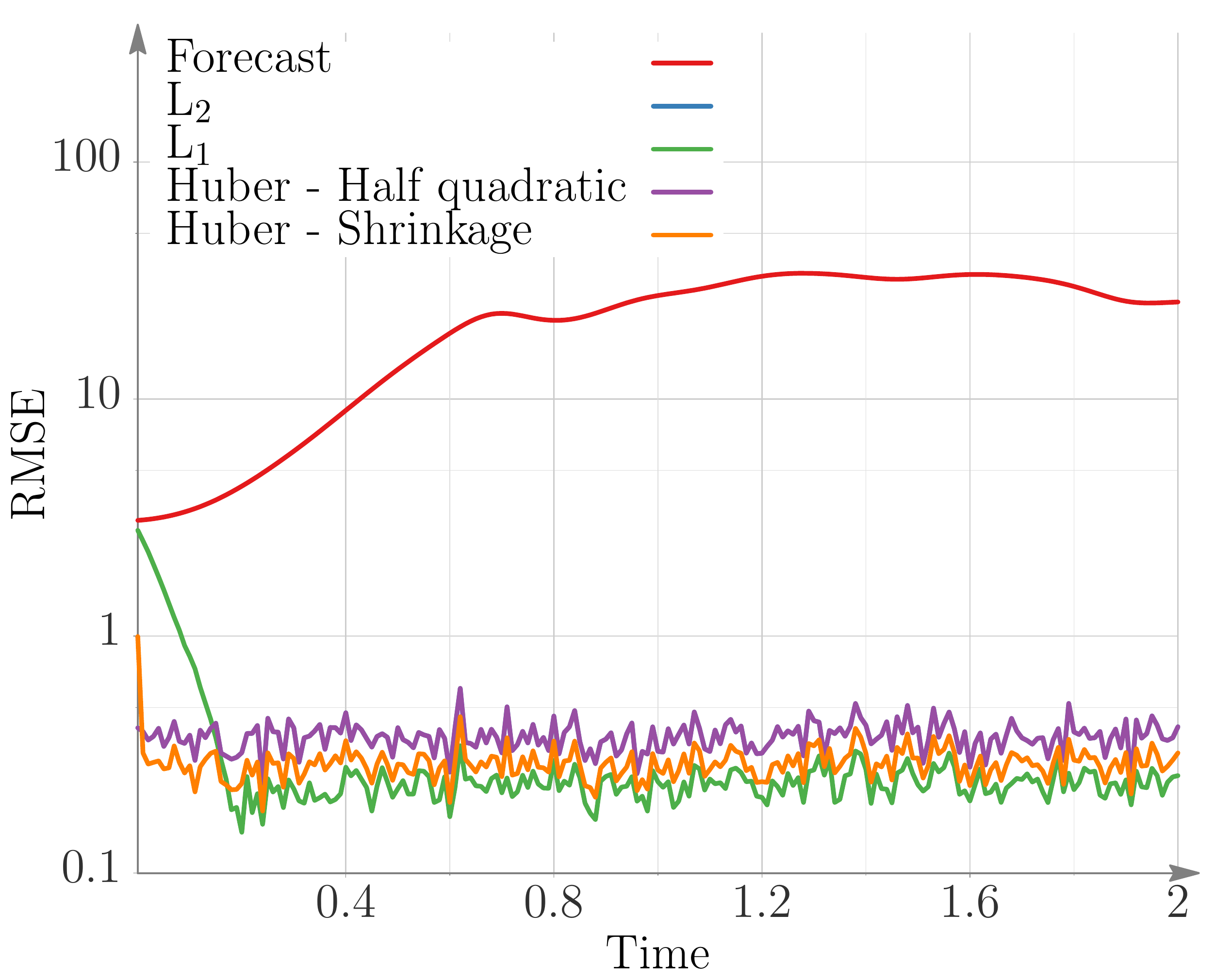}}
  \subfigure[Observations with outliers. Observations taken every 0.01 time units. \label{fig:3DVarRMSE_Lorenz_Bad_0_01_3}]
  {\includegraphics[width=\myfigurewidth]{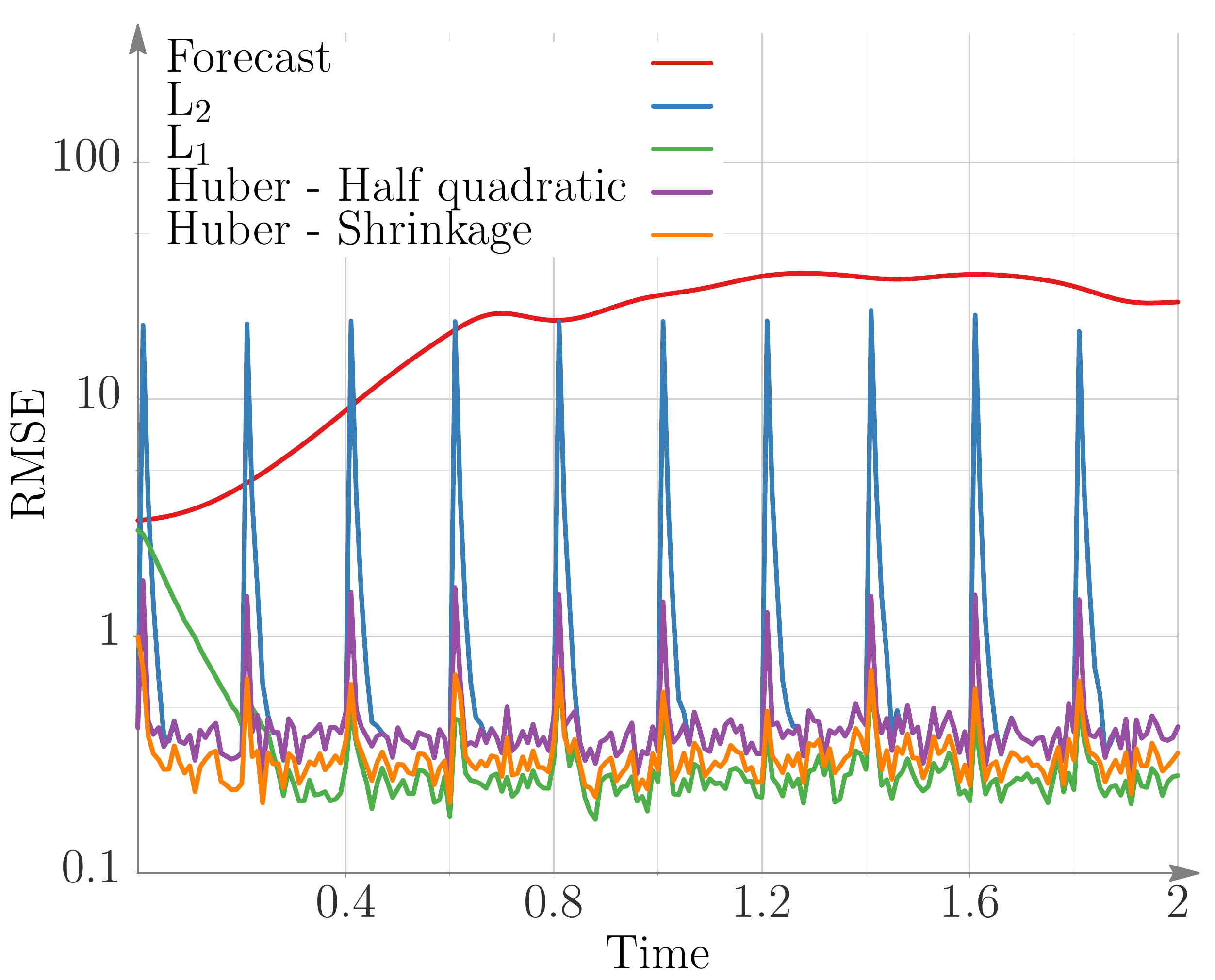}}
  \caption{3D-Var results for the Lorenz-96 model \eqref{eqn:Lorenz96}. The frequency of observations is 0.01 time units. Erroneous observations occur every 0.2 time units. The Huber norm uses $\tau = 3$.}
  \label{fig:Lorenz3DVar_0_01_3}
\end{figure}

\begin{figure}[ht]
  \centering
  \subfigure[Observations with only small random errors. Observations taken every 0.1 time units. \label{fig:3DVarRMSE_Lorenz_Good_0_1_3}]
  {\includegraphics[width=\myfigurewidth]{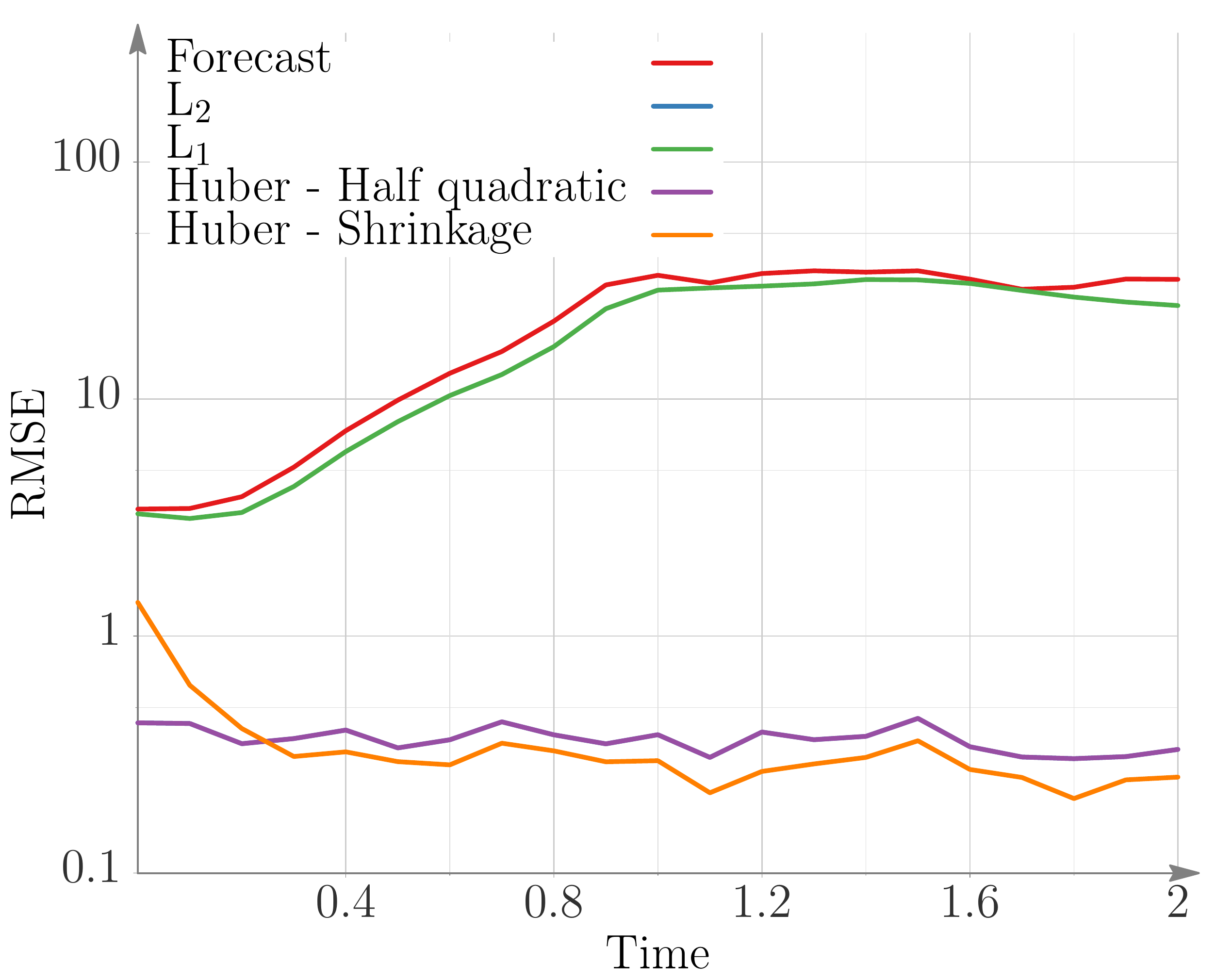}}
  \subfigure[Observations with outliers. Observations taken every 0.1 time units. \label{fig:3DVarRMSE_Lorenz_Bad_0_1_3}]
  {\includegraphics[width=\myfigurewidth]{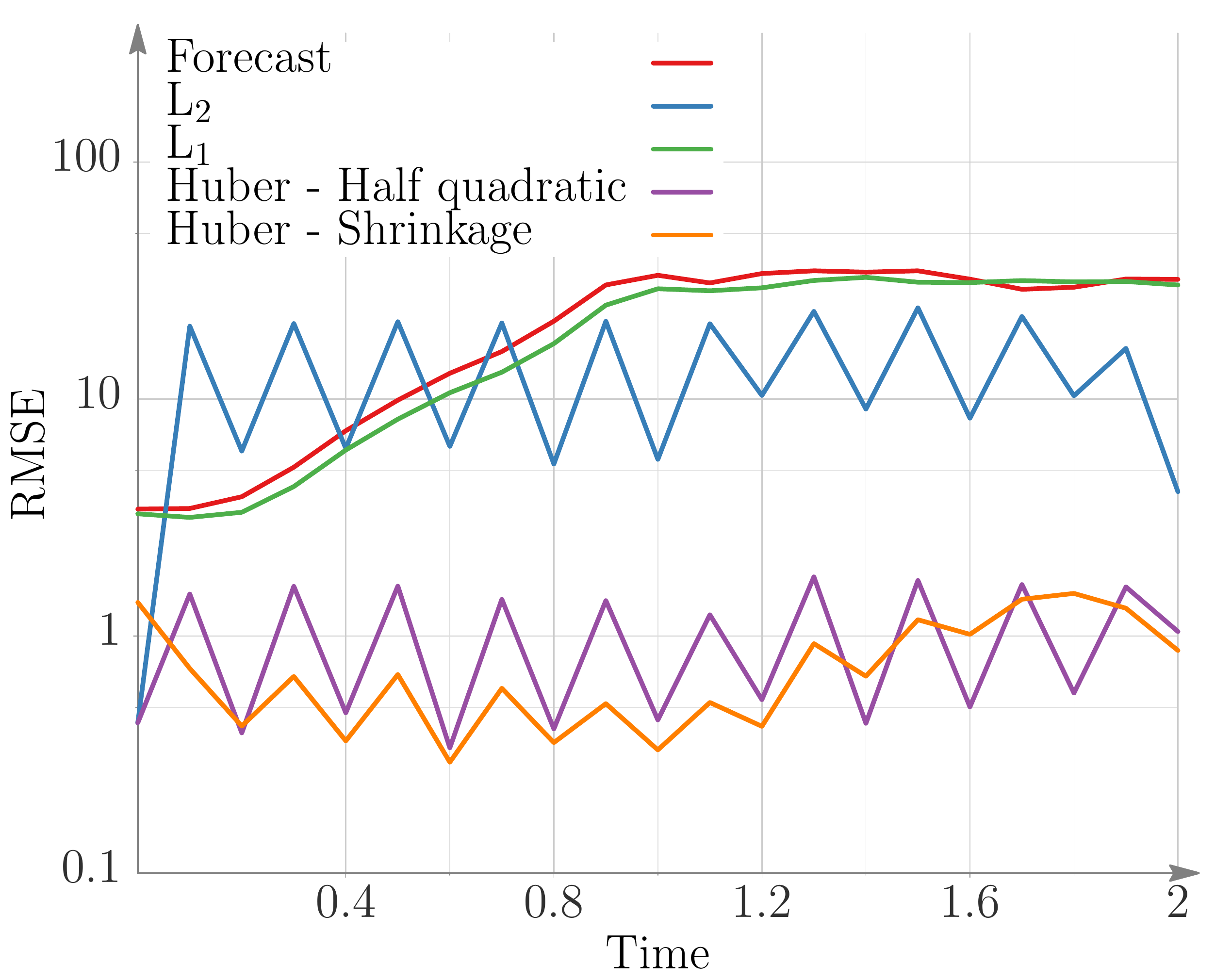}}
  \caption{3D-Var results for the Lorenz-96 model \eqref{eqn:Lorenz96}. The frequency of observations is 0.1 time units. Erroneous observations occur every 0.2 time units. The Huber norm uses $\tau = 3$.}
  \label{fig:Lorenz3DVar_0_1_3}
\end{figure}

The optimization proceeds in stages of inner and outer iterations. The inner iterations aim to find an iterate which minimizes the cost function for a particular value of the penalty parameter ($\mu$). The outer iterations involves increasing the penalty parameter and updating the Lagrange multipliers. For experimental uniformity, we use $15$ outer iterations with 3D-Var, 4D-Var, and LETKF for all formulations: $\Lone$, Huber with ADMM, and Huber with half-quadratic solutions.
Figures \ref{fig:3DVarRMSE_Lorenz_Good_0_01_1} and \ref{fig:3DVarRMSE_Lorenz_Good_0_01_3} show that all formulations perform very well at high observation frequency (one observation every 0.01 units) and using good data. At a low observation frequency with good data the $\Lone$ formulation does not perform as well as the two Huber formulations, as seen in Figures \ref{fig:3DVarRMSE_Lorenz_Good_0_1_1} and \ref{fig:3DVarRMSE_Lorenz_Good_0_1_3}.

For frequent (good) observations both $\Lone$ and $\Ltwo$ formulations perform well, even when outliers are present in the data, as seen in Figures \ref{fig:3DVarRMSE_Lorenz_Bad_0_01_1} and \ref{fig:3DVarRMSE_Lorenz_Bad_0_01_3}. However, at low frequency of measurements and when outliers are present the $\Ltwo$ formulation fails to produce good results, while the Huber formulations demonstrate their robustness, as seen in Figures \ref{fig:3DVarRMSE_Lorenz_Bad_0_1_1} and \ref{fig:3DVarRMSE_Lorenz_Bad_0_1_3}.

\subsection{3D-Var experiments with the shallow water on the sphere model}

Experiments with the shallow water model are performed for an assimilation window of six hours, with all variables being observed hourly. The good data contains random noise but no outliers. The erroneous data contains one outlier, namely the height component at latitude 32$^\circ$ South, longitude 120$^\circ$ East and it has a value which is $\sim$ 50 standard deviations away from the mean. This outlier occurs at all observation times and hence simulates one faulty sensor periodically providing incorrect observation values. 

\begin{figure}[ht]
  \centering
  \subfigure[Observations with small random errors \label{fig:3DVarRMSE_SWE_Good}]
  {\includegraphics[width=\myfigurewidth,height=\myfigureheight]{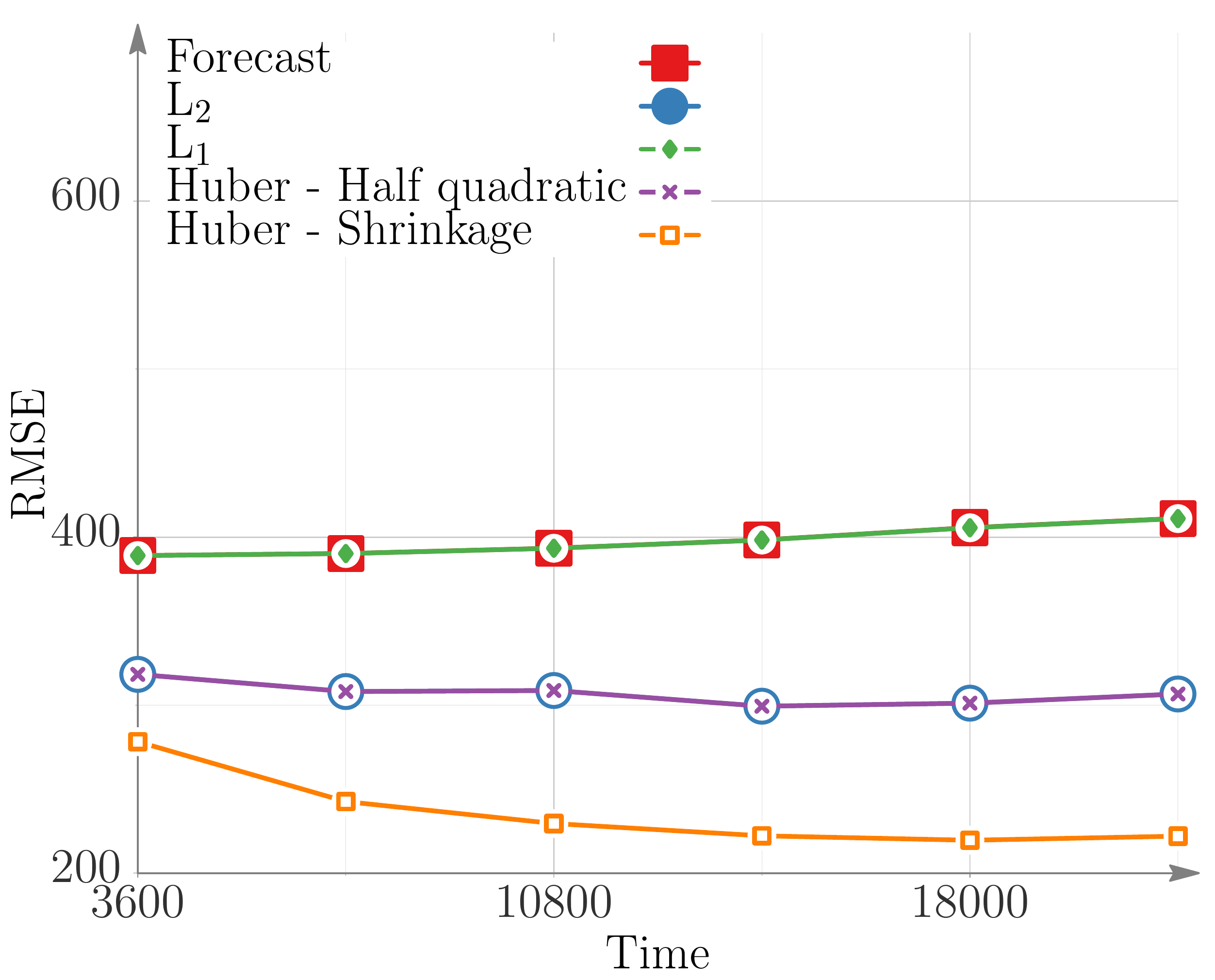}}
  \subfigure[Observations with outliers \label{fig:3DVarRMSE_SWE_Bad}]
  {\includegraphics[width=\myfigurewidth,height=\myfigureheight]{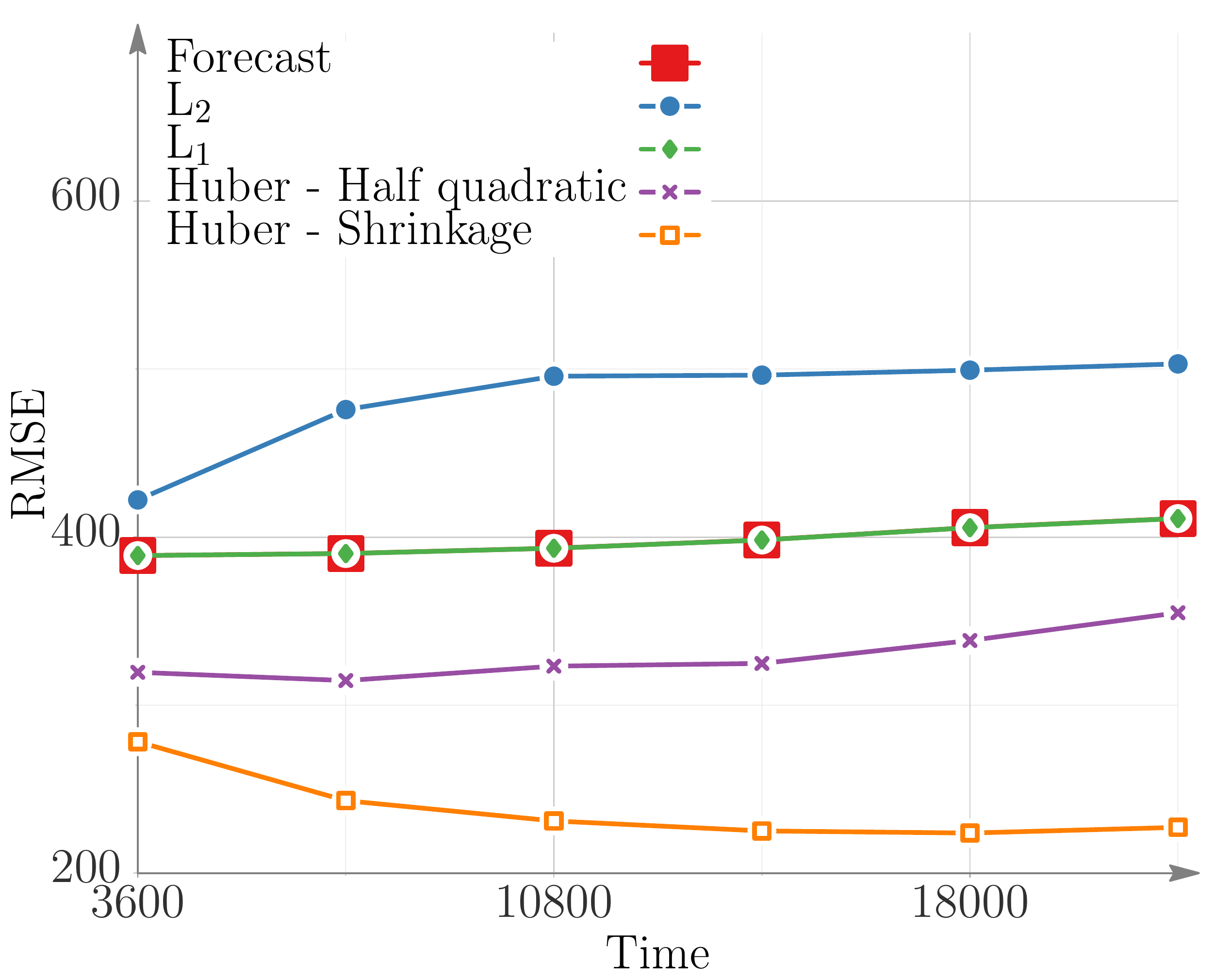}}
  \caption{3D-Var results for the shallow water model \eqref{eqn:swe}. The Huber norm uses $\tau = 2$.}
  \label{fig:SWE3DVar}
\end{figure}

Figures \ref{fig:3DVarRMSE_SWE_Good} and \ref{fig:3DVarRMSE_SWE_Bad} present the results of assimilating good and erroneous observations ,respectively. Root mean square errors (RMSE) are shown for trajectories corresponding to forecast, $\Lone$, $\Ltwo$, and two implementations of Huber norm. The Huber norm implementations perform as well as the $\Ltwo$ norm analysis when only good observations are used. When the observations contain outliers, the $\Ltwo$ norm analysis is inaccurate. However, the Huber norm analysis is robust and remains is unaffected by the quality of observations. The $\Lone$ norm analysis also remains unaffected by the quality of observations, however, its errors are similar to those of the forecast; the $\Lone$ formulation does not use effectively the information from the good observations. Consequently, the Huber norm offers the best 3D-Var formulation, leading to analyses that are both accurate and robust. 
%
%
%

%


\subsection{4D-Var experiments with the Lorenz-96 model}
4D-Var experiments with Lorenz model are performed for one assimilation window of 0.6 units. All components are observed every 0.1 units. The good data contains random noise but no outliers. The erroneous data contains one outlier at all observation times and its value $\sim$ 100 standard deviations away from the mean. 

\begin{figure}[ht]
  \centering
  \subfigure[Observations with small random errors \label{fig:4DVarRMSE_Lorenz_Good}]
  {\includegraphics[width=\myfigurewidth,height=\myfigureheight]{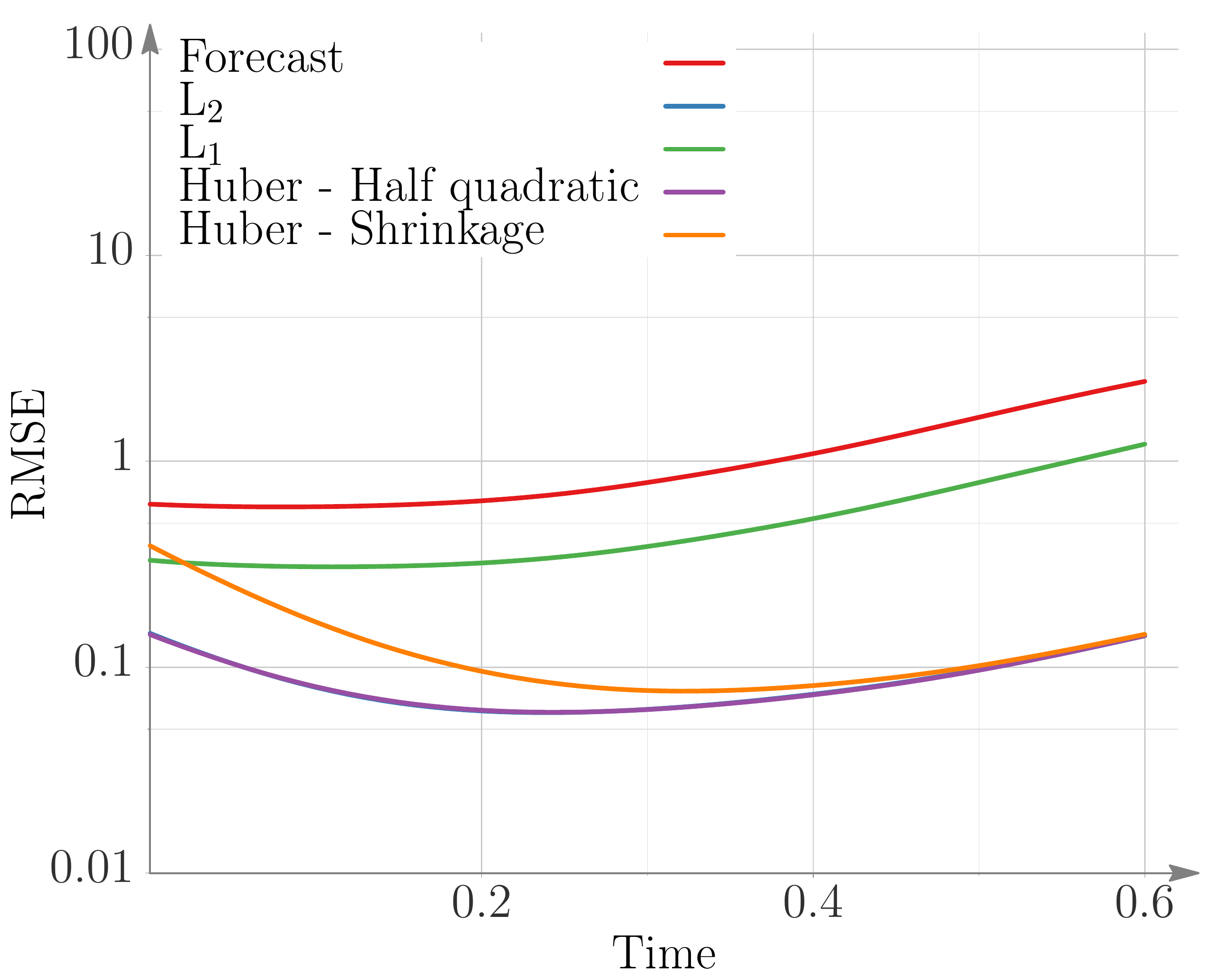}}
  \subfigure[Observations with outliers \label{fig:4DVarRMSE_Lorenz_Bad}]
  {\includegraphics[width=\myfigurewidth,height=\myfigureheight]{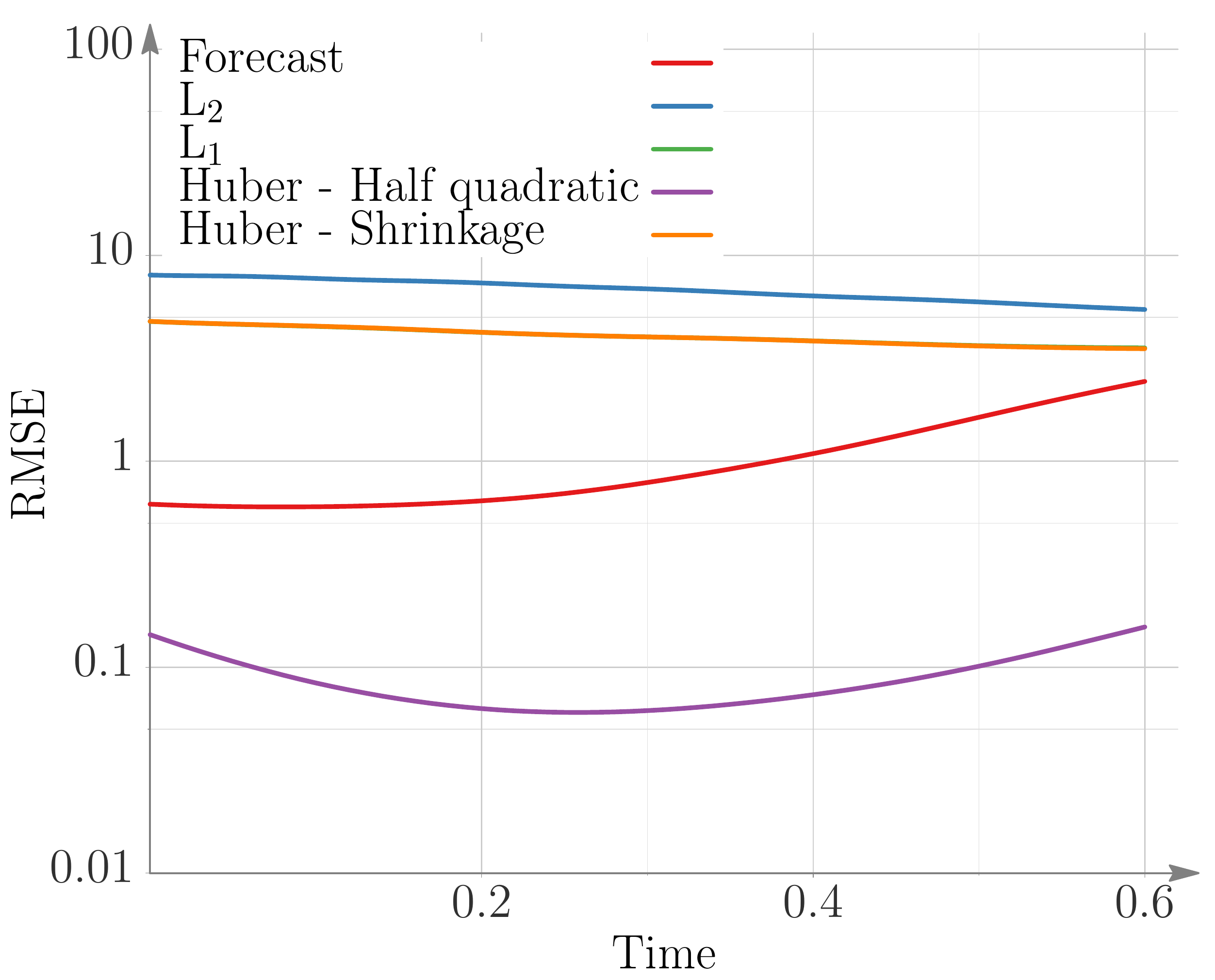}}
  \caption{4D-Var results for the Lorenz-96 model \eqref{eqn:Lorenz96}. The Huber norm uses $\tau = 2$.}
  \label{fig:Lorenz4DVar}
\end{figure}

Figure \ref{fig:4DVarRMSE_Lorenz_Good} presents the 4D-Var results with good data. The performance of the Huber-4D-Var with half-quadratic formulation matches the performance of $\Ltwo$ formulation. The analyses provided by $\Lone$ and Huber formulations with shrinkage are less accurate. Figure \ref{fig:4DVarRMSE_Lorenz_Bad} presents the results of assimilating data with outliers. The 4D-Var-Huber analysis using the half quadratic formulation remains unaffected by the data errors. All other formulations ($\Lone$, $\Ltwo$, and Huber with shrinkage) are negatively impacted by the presence of outliers.

\subsection{4D-Var experiments with the shallow water on the sphere model}

We consider one nine hours long assimilation window, with all components of the system observed hourly.
The good data contains only small random errors. The bad data contains one large outlier with an error $\sim50$ standard deviations away from the mean. This simulates a faulty sensor that provides erroneous data at all observation times.

\begin{figure}[ht]
  \centering
  \subfigure[Observations with small random errors \label{fig:4DVarRMSE_SWE_Good}]
  {\includegraphics[width=\myfigurewidth,height=\myfigureheight]{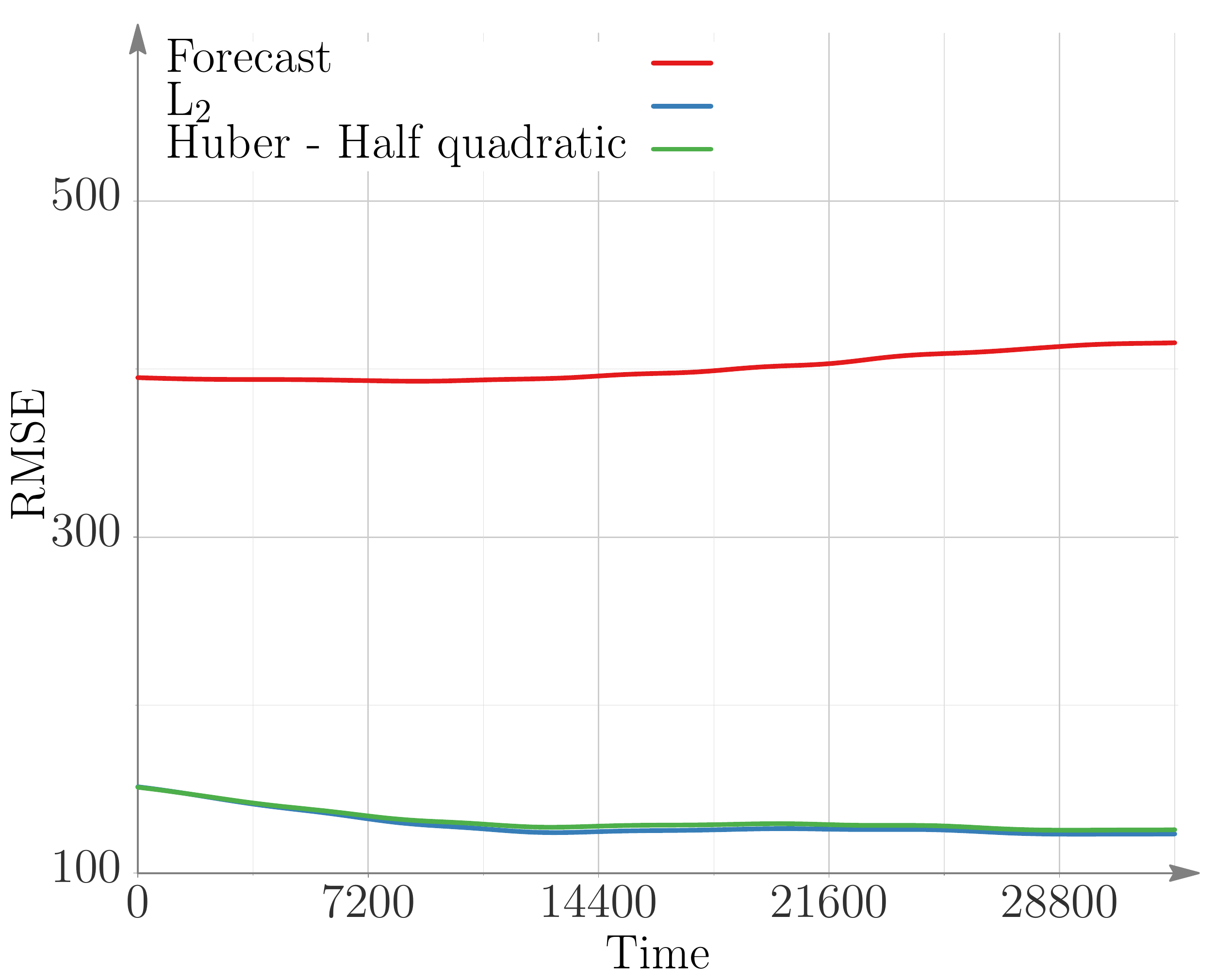}}
  \subfigure[Observations with outliers \label{fig:4DVarRMSE_SWE_Bad}]
  {\includegraphics[width=\myfigurewidth,height=\myfigureheight]{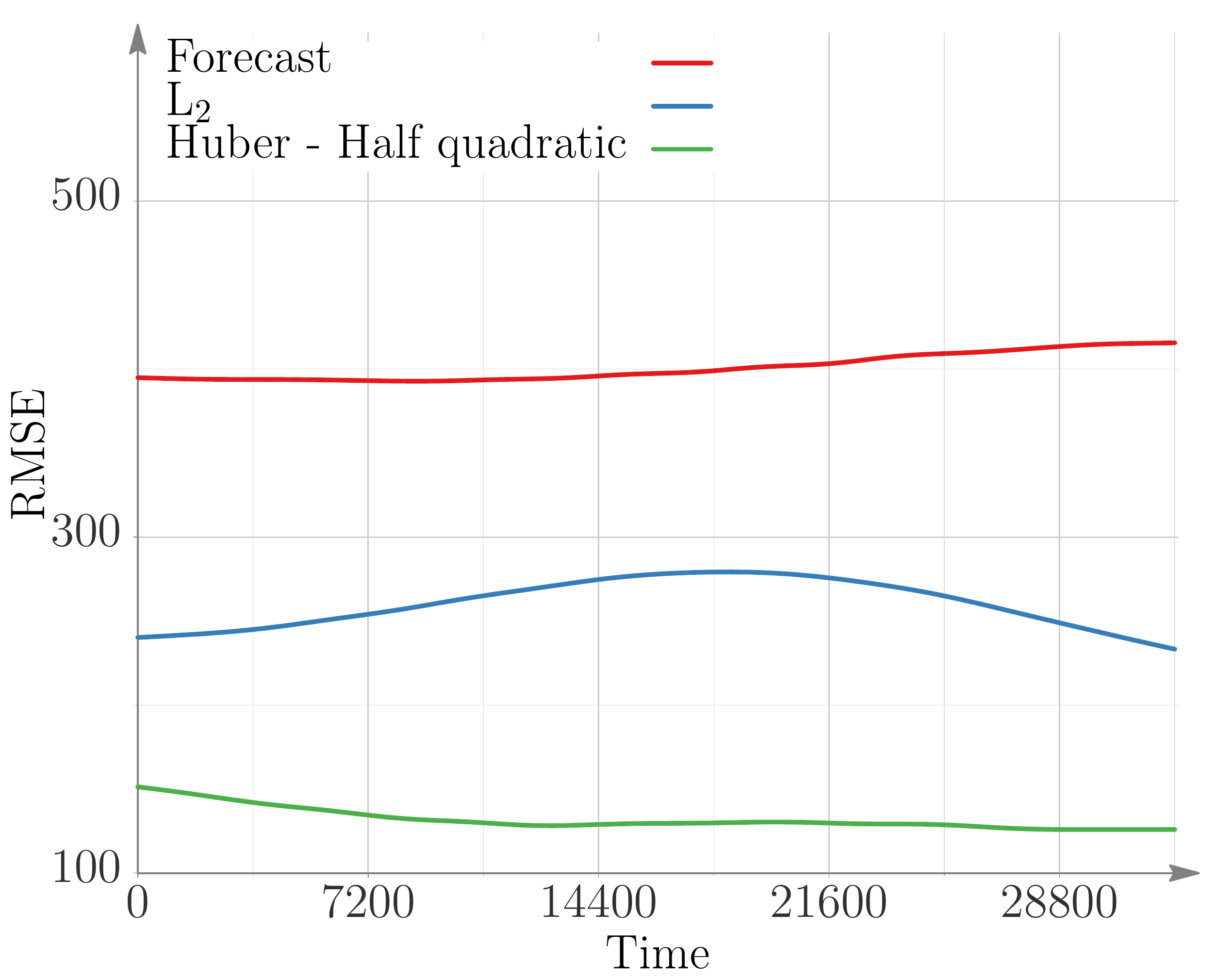}}
  \caption{4D-Var results for the shallow water model \eqref{eqn:swe}. The Huber norm uses $\tau = 2$.}
  \label{fig:SWE4DVar}
\end{figure}

The $\Lone$ and Huber norm with shrinkage formulations of 4D-Var did not perform well for the Lorenz model, and we do not test them on the shallow water model. Figures \ref{fig:4DVarRMSE_SWE_Good} and \ref{fig:4DVarRMSE_SWE_Bad} show the 4D-Var results for the shallow water model using good and faulty observations respectively. We compare the performance of Huber with half-quadratic formulation with the $\Ltwo$ norm formulation.  When the quality of the observations is good the Huber norm and the $\Ltwo$ norm formulations perform equally well. When the observations contain outliers the $\Ltwo$ formulation gives inaccurate analyses, while the Huber norm formulation is robust and continues to provide good results.


\subsection{EnSRF experiments with the Lorenz-96 model}

In this section we discuss the results of the LETKF implementations based on the $\Lone$, $\Ltwo$, and Huber norms. We note that, due to localization the data outliers in LETFK will impact the states in only small regions of the domain. The number of ensemble members for the experiments is 20.

\begin{figure}[ht]
  \centering
  \subfigure[Observations with only small random errors \label{fig:EnKFRMSE_Lorenz_Good_0_01_1}]
  {\includegraphics[width=\myfigurewidth,height=\myfigureheight]{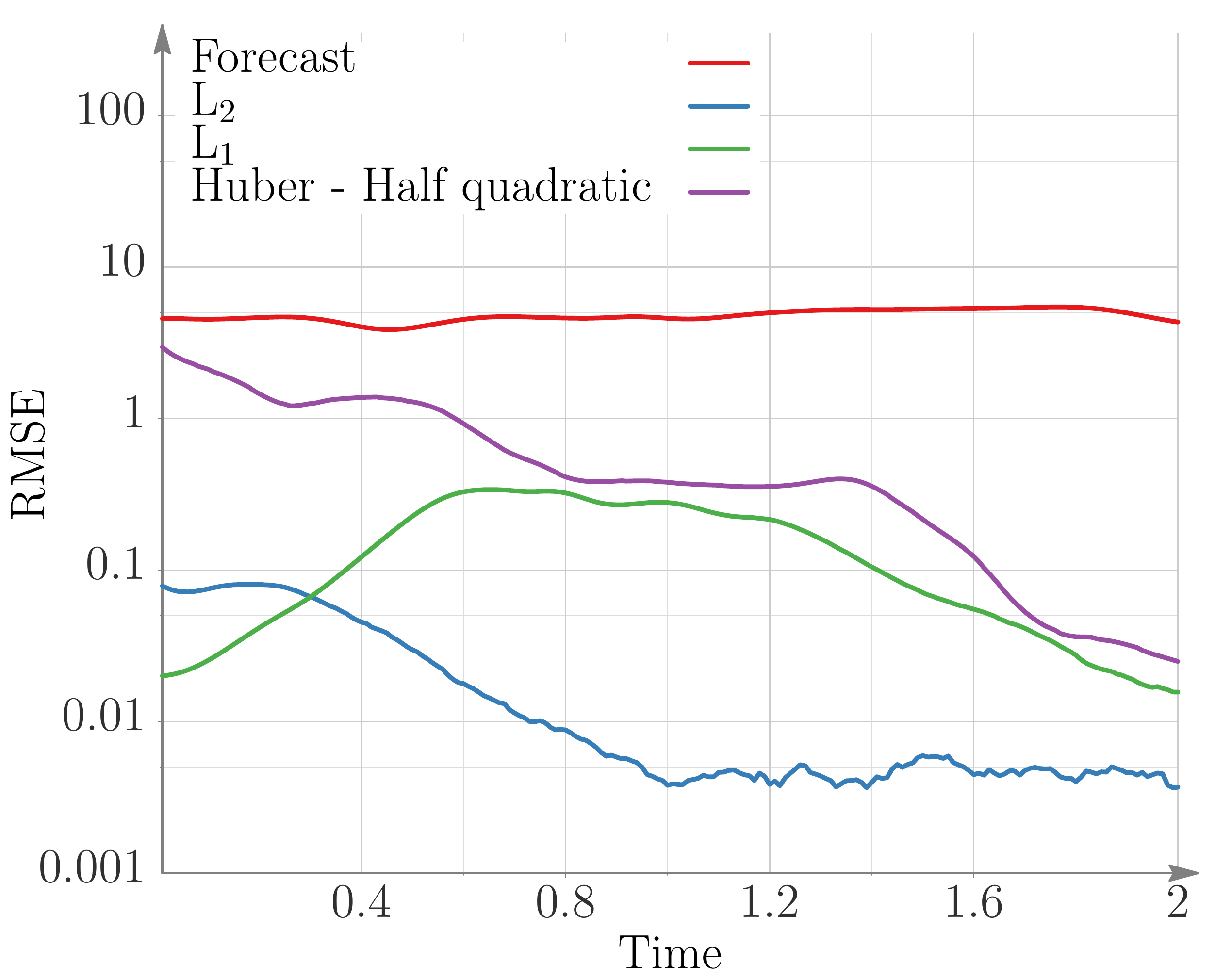}}
  \subfigure[Observations with outliers \label{fig:EnKFRMSE_Lorenz_Bad_0_01_1}]
  {\includegraphics[width=\myfigurewidth,height=\myfigureheight]{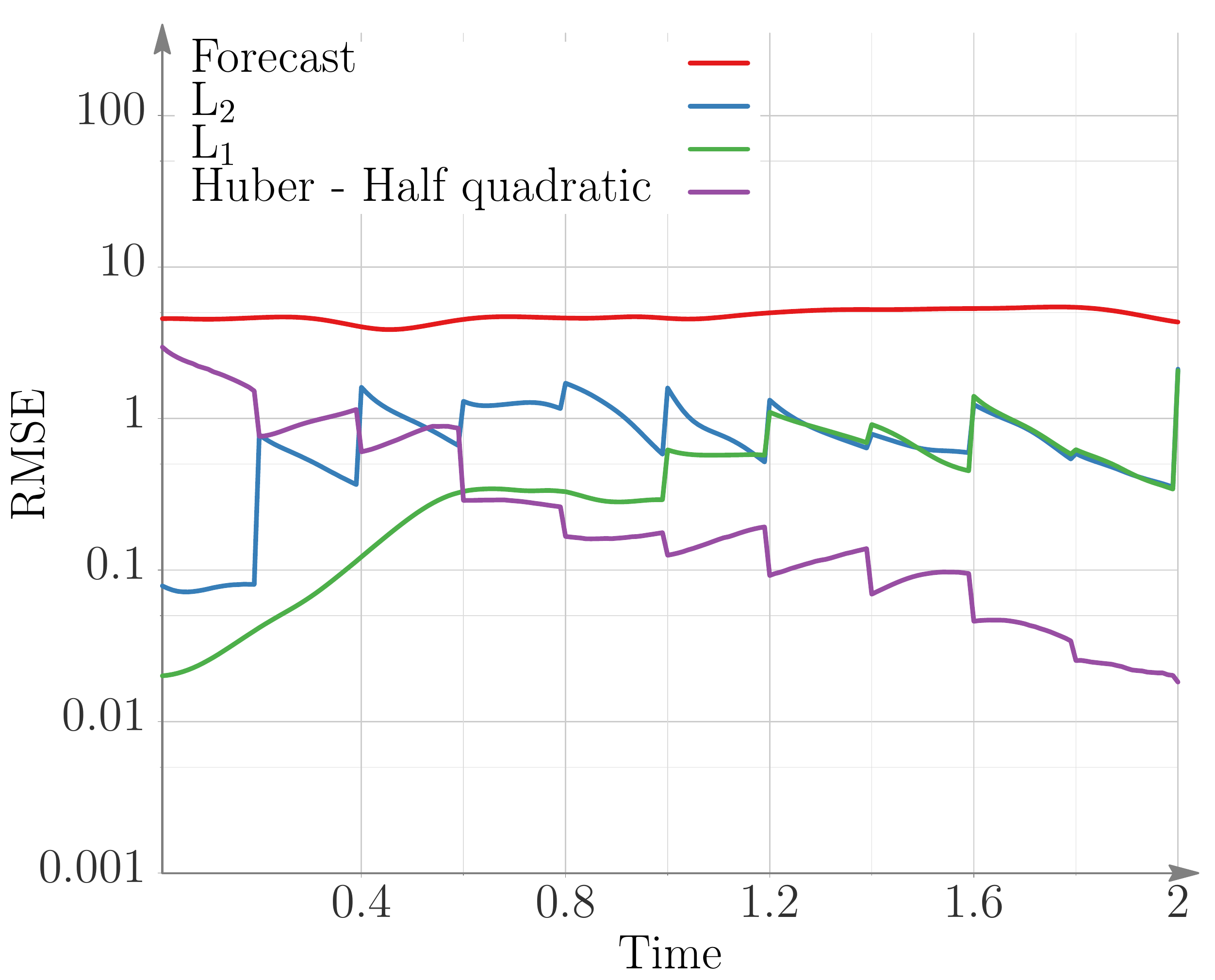}}
  \caption{LETKF results for the Lorenz-96 model \eqref{eqn:Lorenz96}. The frequency of observations is 0.01 time units. Erroneous observations occur every 0.2 time units. The Huber norm uses $\tau = 1$.}
  \label{fig:LorenzEnKF_0_01_1}
\end{figure}

Assimilation experiments are carried out with different frequencies for outliers in the observations. Figures \ref{fig:EnKFRMSE_Lorenz_Good_0_01_1} and \ref{fig:EnKFRMSE_Lorenz_Bad_0_01_1} show results with observations taken every 0.01 time units; for the bad data outliers are present every 0.2 time units. The $\Ltwo$ formulation converges the fastest when only good data is used, as seen in Figure \ref{fig:EnKFRMSE_Lorenz_Good_0_01_1}. The $\Lone$ and Huber formulations converge slower, and give similar results at the end of the assimilation window.  When bad data is used LETKF-Huber outperforms the other formulations, as illustrated in Figure \ref{fig:EnKFRMSE_Lorenz_Bad_0_01_1}. The $\Ltwo$ LETKF formulation still gives improvements over the forecast, since the frequency of outliers is small (every 20 observation times), and injection of good observations at other times avoids filter divergence. 

\begin{figure}[ht]
  \centering
  \subfigure[Observations with only small random errors \label{fig:EnKFRMSE_Lorenz_Good_0_1_1}]
  {\includegraphics[width=\myfigurewidth,height=\myfigureheight]{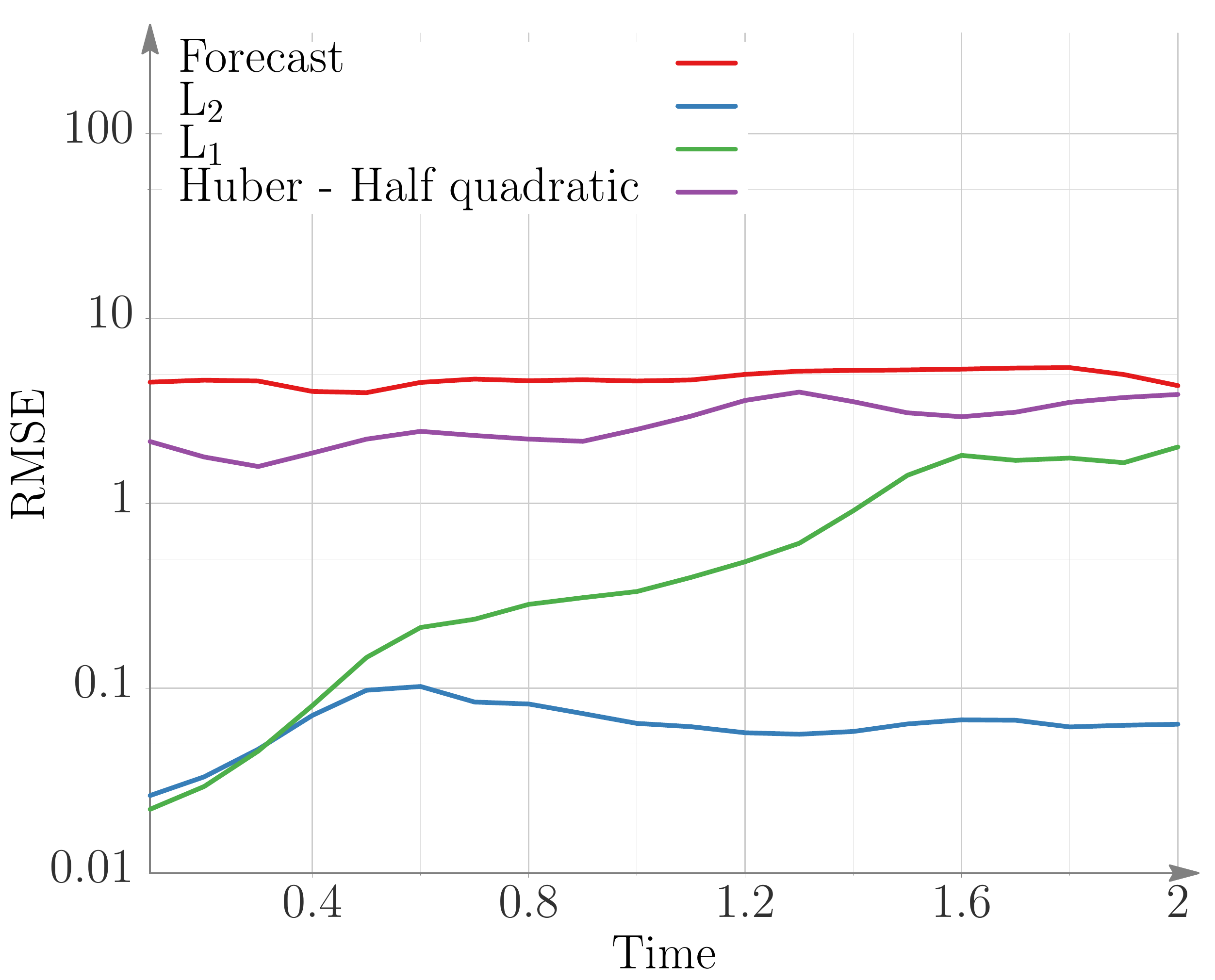}}
  \subfigure[Observations with outliers \label{fig:EnKFRMSE_Lorenz_Bad_0_1_1}]
  {\includegraphics[width=\myfigurewidth,height=\myfigureheight]{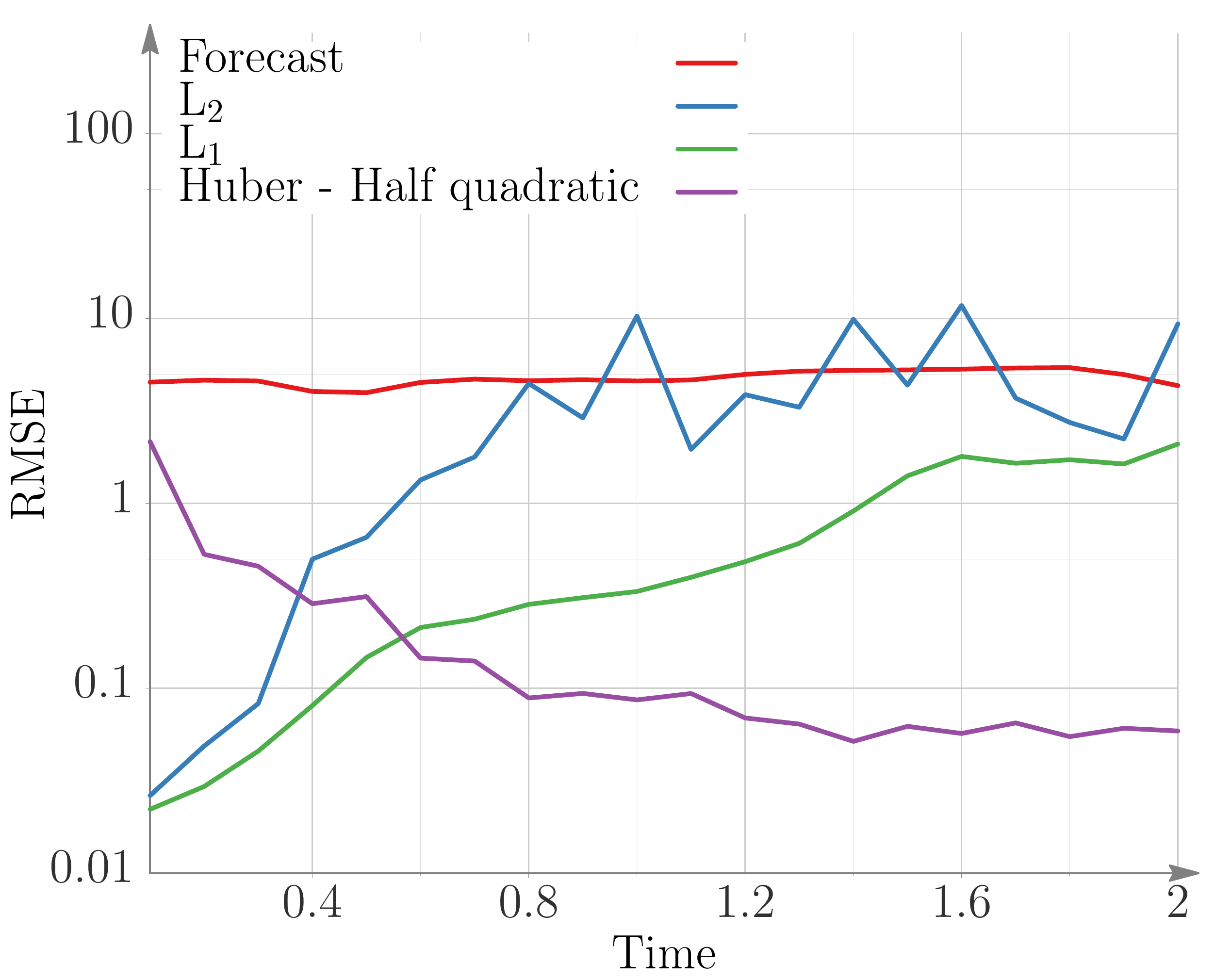}}
  \caption{LETKF results for the Lorenz-96 model \eqref{eqn:Lorenz96}. The frequency of observations is 0.1 time units. Erroneous observations occur every 0.2 time units. The Huber norm uses $\tau = 1$.}
  \label{fig:LorenzEnKF_0_1_1}
\end{figure}

\begin{figure}[ht]
  \centering
  \subfigure[Observations with only small random errors \label{fig:EnKFRMSE_Lorenz_Good_0_01_3}]
  {\includegraphics[width=\myfigurewidth,height=\myfigureheight]{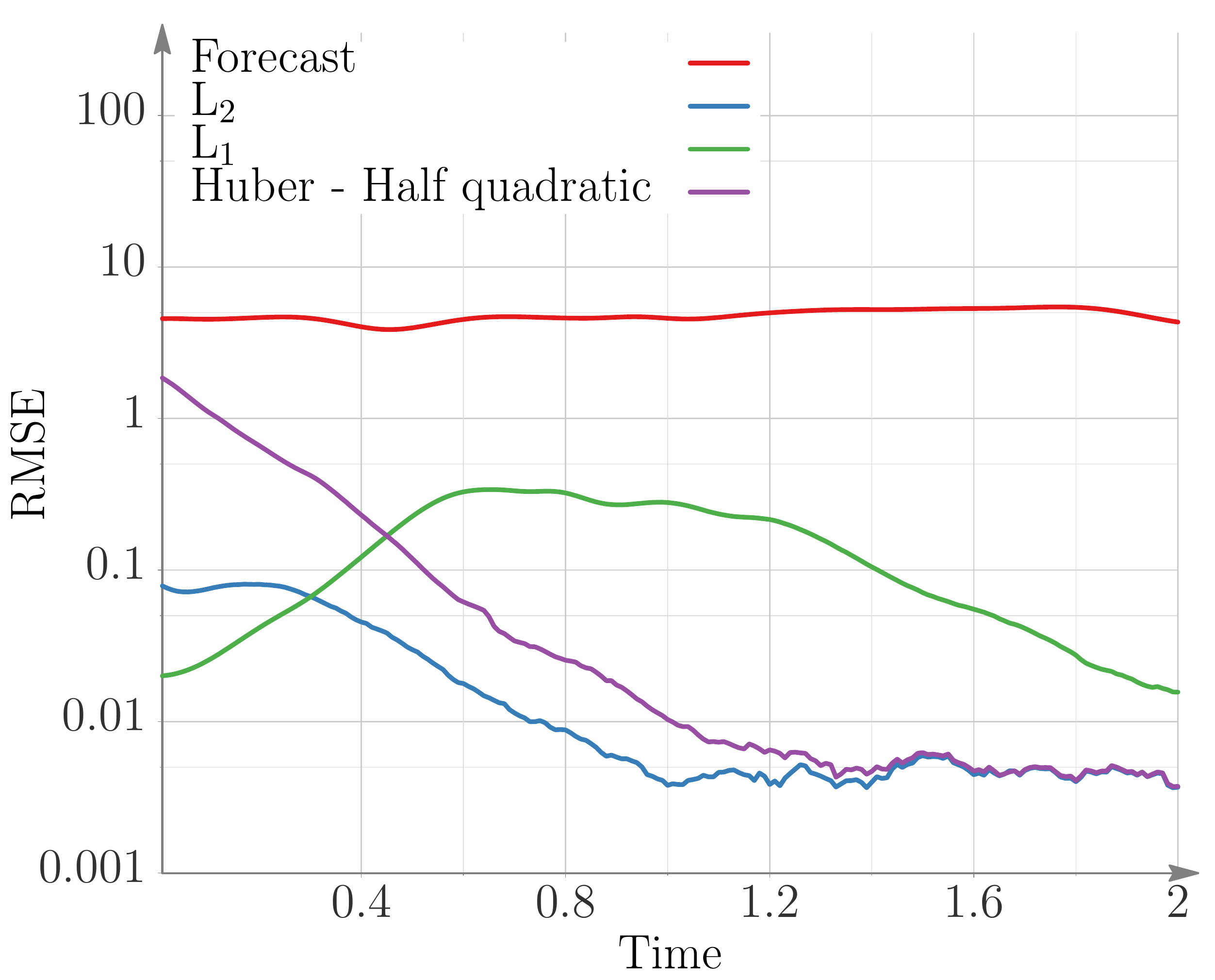}}
  \subfigure[Observations with outliers \label{fig:EnKFRMSE_Lorenz_Bad_0_01_3}]
  {\includegraphics[width=\myfigurewidth,height=\myfigureheight]{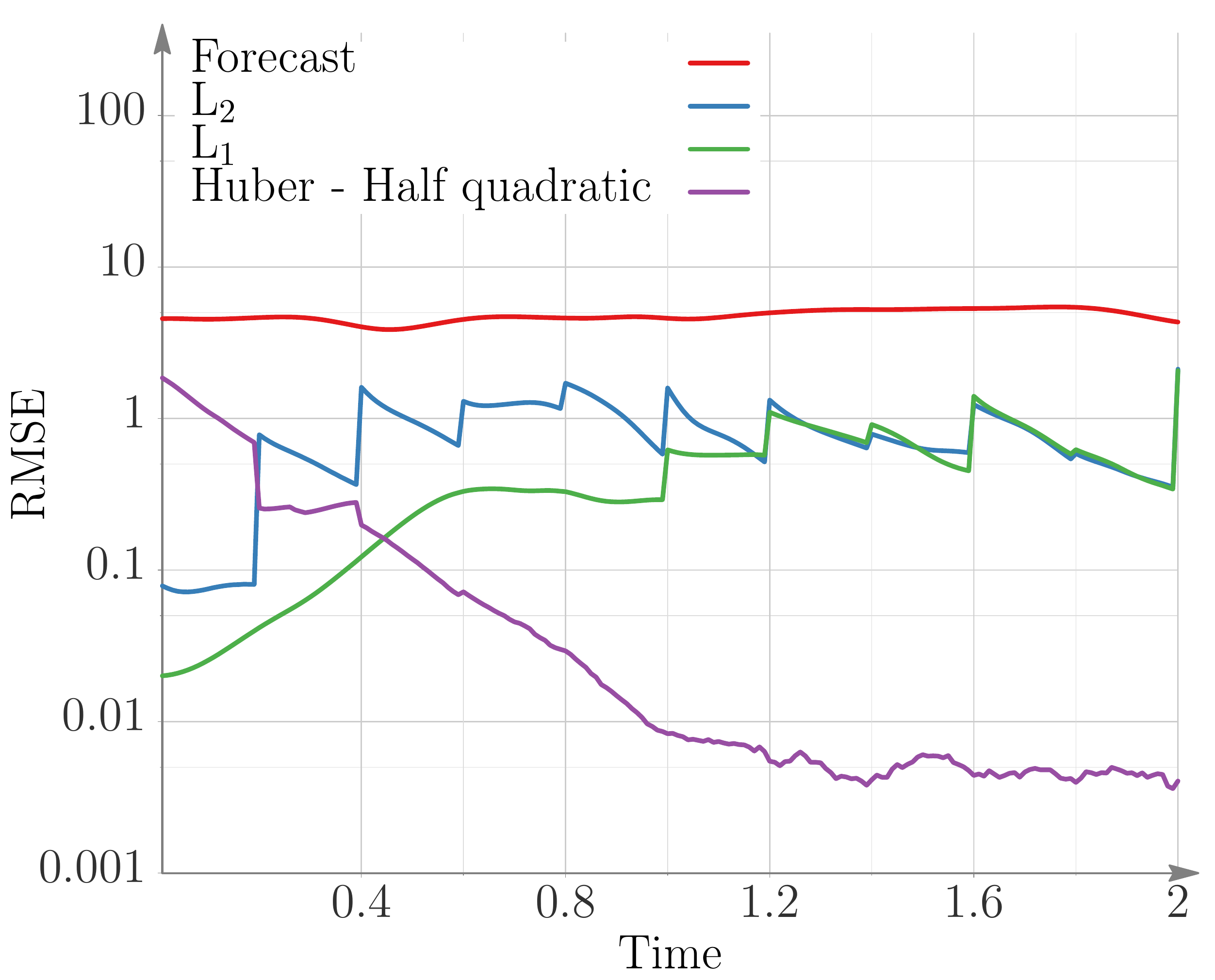}}
  \caption{LETKF results for the Lorenz-96 model \eqref{eqn:Lorenz96}. The frequency of observations is 0.01 time units. Erroneous observations occur every 0.2 time units. The Huber norm uses $\tau = 3$.}
  \label{fig:LorenzEnKF_0_01_3}
\end{figure}

Figures \ref{fig:EnKFRMSE_Lorenz_Good_0_1_1} and \ref{fig:EnKFRMSE_Lorenz_Bad_0_1_1} show results with observations taken every 0.1 time units; for the bad data outliers are present every 0.2 time units. The relatively frequent outliers (every other observation) affects the analysis obtained with the $\Ltwo$ EnKF formulation. In this case the $\Lone$ formulation performs better that the $\Ltwo$ one. Huber LETKF gives the best analyses overall. When $\tau=1$, the method can treat good observations as outliers, and therefore the contributions of some good observations to the analyses can be attenuated; this can be seen in Figures \ref{fig:EnKFRMSE_Lorenz_Good_0_1_1} and \ref{fig:EnKFRMSE_Lorenz_Bad_0_1_1} where the Huber formulation is less accurate than the $\Ltwo$ formulation for good data. When $\tau$ is increased, for instance $\tau=3$, the trajectories of the Huber LETKF for good and bad observations are almost identical, as can be seen in Figures \ref{fig:EnKFRMSE_Lorenz_Good_0_01_3} and \ref{fig:EnKFRMSE_Lorenz_Bad_0_01_3}.

\begin{figure}[ht]
  \centering
  \subfigure[Observations with only small random errors \label{fig:EnKFRMSE_Lorenz_Good_0_1_3}]
  {\includegraphics[width=\myfigurewidth,height=\myfigureheight]{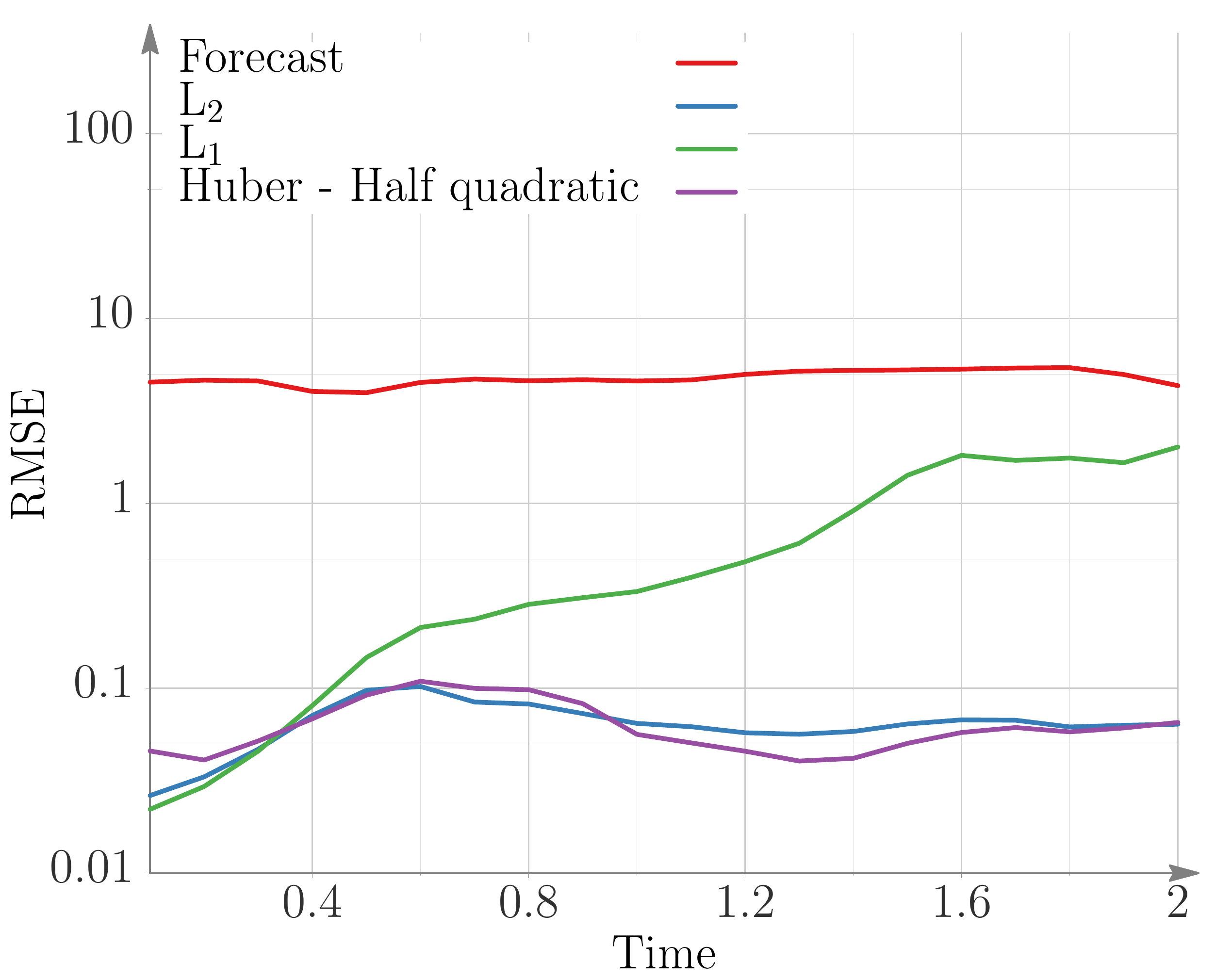}}
  \subfigure[Observations with outliers \label{fig:EnKFRMSE_Lorenz_Bad_0_1_3}]
  {\includegraphics[width=\myfigurewidth,height=\myfigureheight]{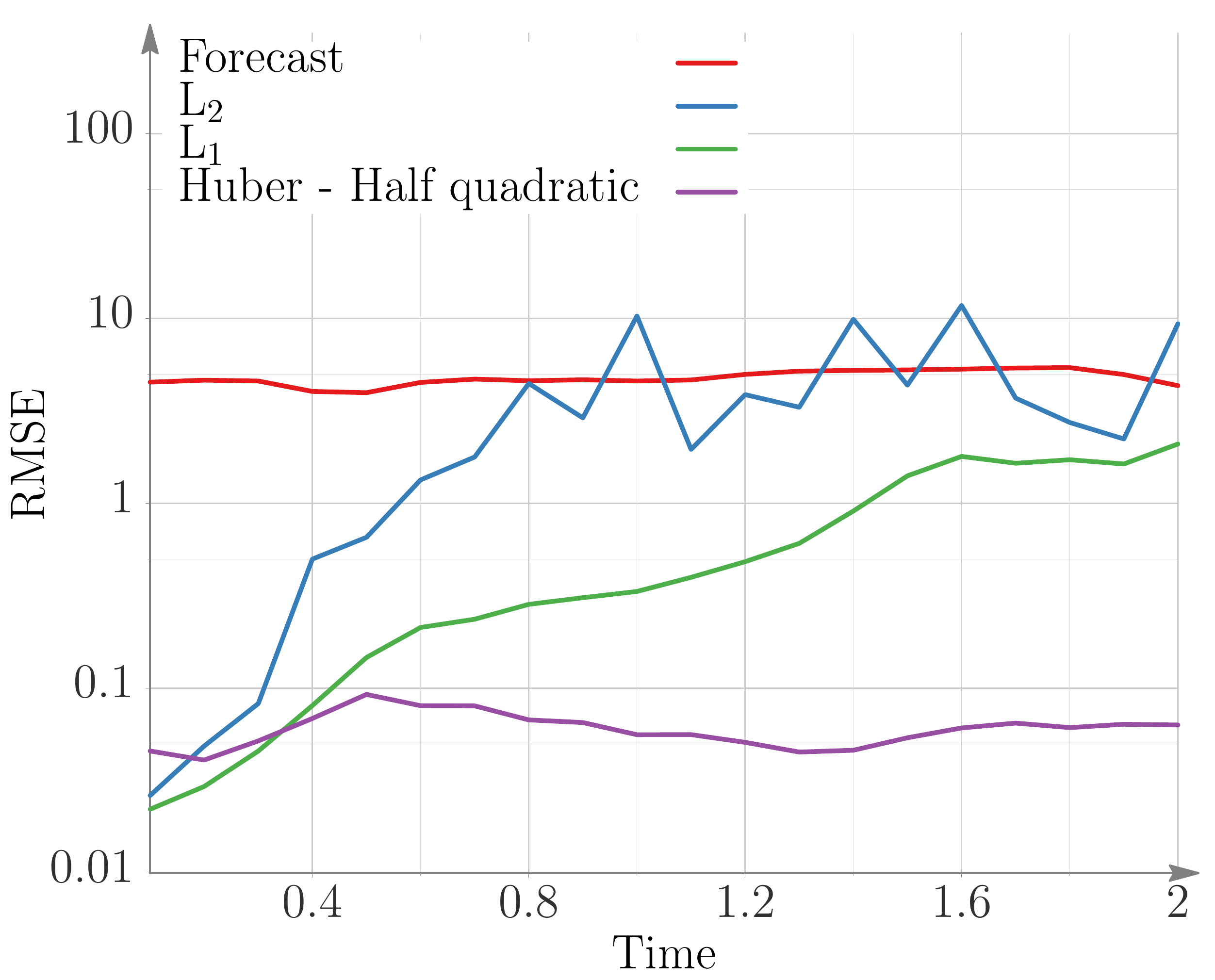}}
  \caption{LETKF results for the Lorenz-96 model \eqref{eqn:Lorenz96}. The frequency of observations is 0.1 time units. Erroneous observations occur every 0.2 time units. The Huber norm uses $\tau = 3$.}
  \label{fig:LorenzEnKF_0_1_3}
\end{figure}


\subsection{EnSRF experiments with the shallow water on the sphere model}

\begin{figure}[ht]
  \centering
  \subfigure[Observations with small random errors \label{fig:EnKFRMSE_SWE_Good}]
  {\includegraphics[width=\myfigurewidth,height=\myfigureheight]{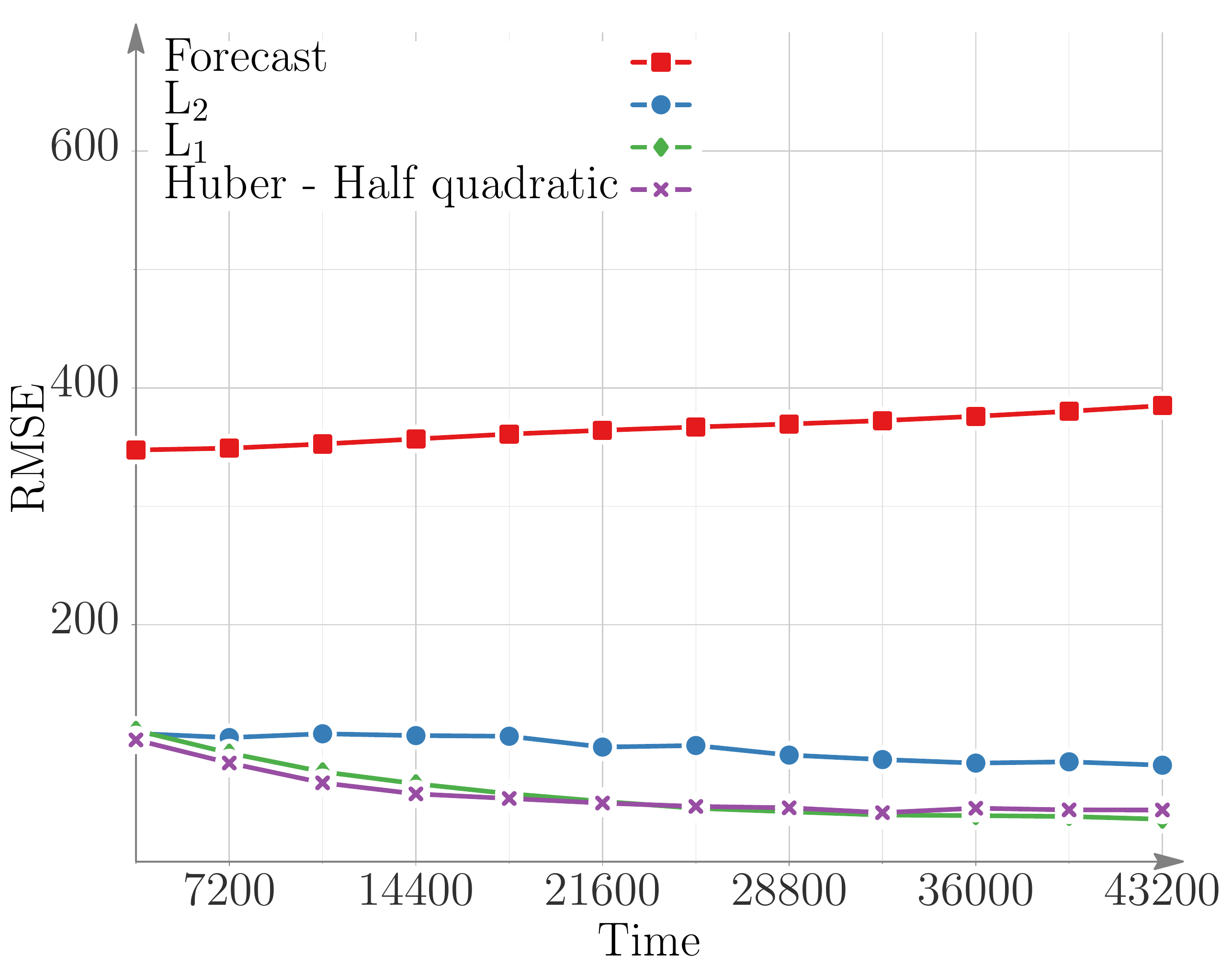} \label{ed:no-outliers}}
  \subfigure[Observations with outliers \label{fig:EnKFRMSE_SWE_Bad}]
  {\includegraphics[width=\myfigurewidth,height=\myfigureheight]{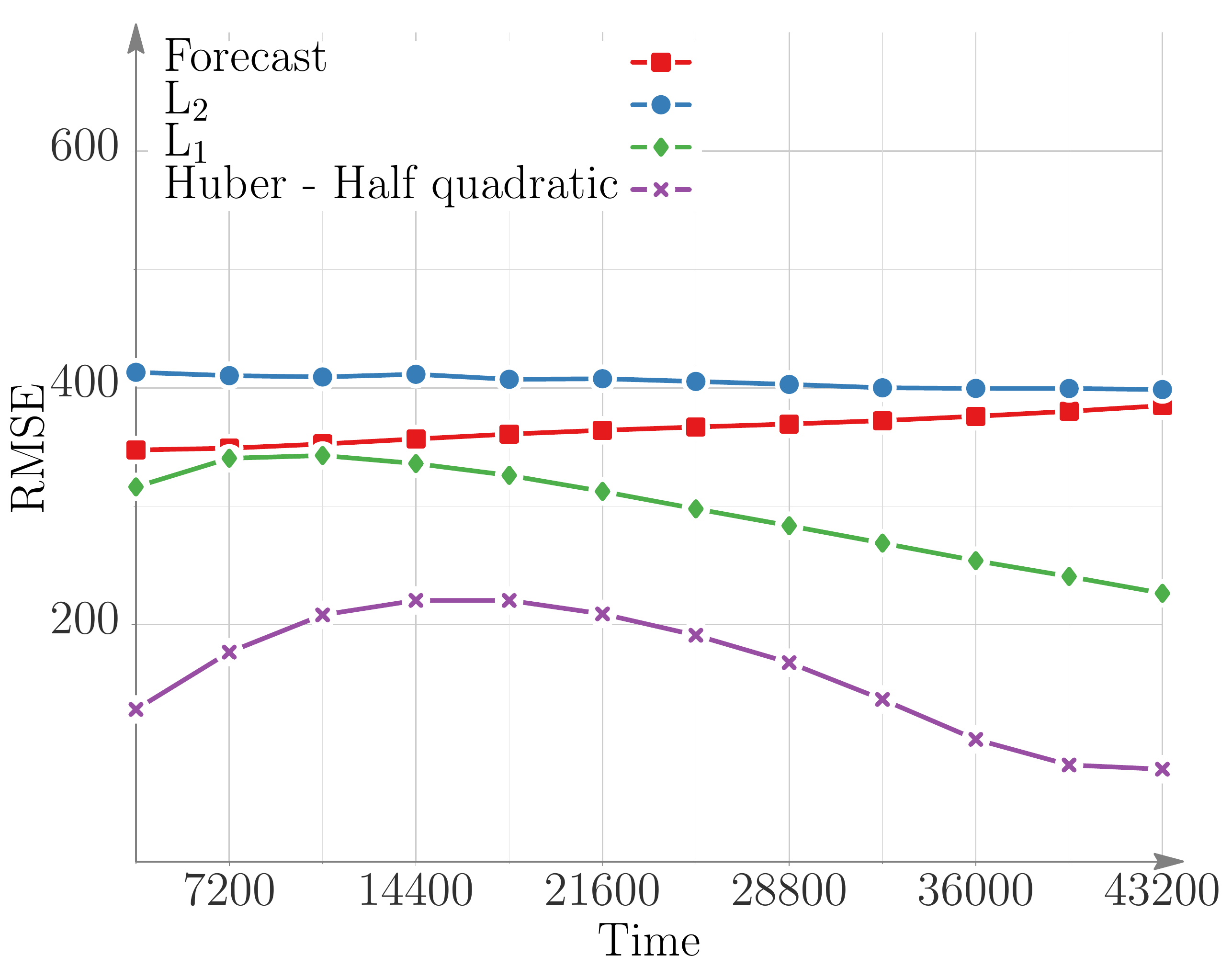} \label{ed:outliers}}
\caption{LETKF results for the shallow water model \eqref{eqn:swe}. The Huber norm uses $\tau = 1$.}
  \label{fig:SWE-Model-EnKF}
\end{figure}

The shallow water model \eqref{eqn:swe} experiments use ensembles with 30 members. The initial ensemble perturbation is normal, with a standard deviation of 5\% of the background state for each component. Figure \ref{fig:SWE-Model-EnKF} shows the analysis errors for the $\Lone$, $\Ltwo$, and Huber LETKF formulations. The Huber half quadratic LETKF provides analyses that are as accurate as the $\Ltwo$ formulation analyses when only good data is used. When outliers are present, however, the $\Ltwo$ results are very inaccurate. The Huber formulation is almost completely unaffected by the data outliers, and provides equally better analyses with and without data outliers. The $\Lone$ formulation results are also unaffected by outliers, but their overall accuracy is quite low.

\section{Conclusions} \label{sec:conc}
This papers develops a systematic framework for performing robust 3D-Var, 4D-Var, and ensemble-based filtering data assimilation. The traditional algorithms are formulated as optimizations problems where the cost functions are $\Ltwo$ norms of background and observation residuals; the $\Ltwo$ norm choice corresponds to Gaussianity assumptions for background and observation errors, respectively. The $\Ltwo$ norm formulation has to make large corrections to the optimal solution in order to accommodate outliers in the data. Consequently, a few observations containing large errors can considerably deteriorate the overall quality of the analysis. In order to avoid this effect traditional data quality control rejects certain observations, which may result in loss of useful information. A more recent approach performs a one-time, off-line adjustment of the data weights based on the estimated data quality.

The robust data assimilation framework described herein reformulates the underlying optimization problems using $\Lone$ and Huber norms. The resulting robust data assimilation algorithms do not reject any observations. Moreover, the relative importance weights of the data are not decided off-line, but rather are adjusted iteratively for each observation based on their deviation from the mean forecast.
Numerical results show that the  $\Lone$ norm formulation is very robust, being the least affected by presence of outliers. However, the resulting analyses are inaccurate: the  $\Lone$ norm rejects not only outliers, but useful information as well. The Huber norm formulation is able to fully use the information from ``good data'' while remaining robust (rejecting the influence of outliers). We consider two solution methods for Huber norm optimization. The shrinkage operator solution displays slow convergence and can become impractical in real-life problems. The Huber norm solution using the half-quadratic formulation seems to be the most suitable approach for large scale data assimilation applications. It is only slightly more expensive than the traditional  $\Ltwo$ 4D-Var approach, and yields good results both in the absence and in the presence of data outliers. Future work will apply robust data assimilation algorithms using a Huber norm formulation with a half-quadratic solution to real problems using the Weather Research and Forecasting model \cite{WRF}.

\bibliographystyle{plain}
\bibliography{robust_da}

\appendix
\section{Multivariate Laplace distribution of observation errors}

The univariate Laplace distribution models an exponential decay on both sides away from the mean $\mu$; the probability density function is:
\[
\begin{split}
\mathcal{P}(z)
&= (2\,\lambda)^{-1}\, \exp\left( - \left| z - \mu \right|/\lambda \right), \quad
\tt{E}\left[ z \right] = 0, \quad \tt{Var}\left[ \z \right] = 2\,\lambda^2.
\end{split}
\]
Let $\z$ be the scaled innovation vector \eqref{eqn:define-z}, with $\tt{E}\left[ \z \right] = 0$, and $\tt{Cov}\left[ \z \right] = \mathbf{I}$. Assume that each component error has a univariate Laplace distribution, and that the component errors are independent. The observation part of the variational cost function \eqref{eqn:3dvar-L1} is the log-likelihood:
\[
-\log\; \mathcal{P}(\z) \propto \lambda^{-1} \sum_{\ell=1}^n \left| \z_\ell \right|, \quad
\tt{E}\left[ \z \right] = \mu, \quad \tt{Cov}\left[ \z \right] = 2\,\lambda^2\, \mathbf{I}.
\]
Our formulation of the variational cost function \eqref{eqn:3dvar-L1} uses $\lambda=2$, which translates in an assumed variance for each component of $8$.

The multivariate Laplace distribution is \cite{Eltoft_2006_multivariate-Laplace}
\begin{equation}
\label{eqn:laplace-mv}
\begin{split}
\mathcal{P}(\z) 
 &= \frac{2}{(2\,\lambda)^{(n+2)/4}\, \pi^{n/2} }\, \frac{K_{(n-2)/2}\left( \frac{2}{\lambda}\, \left\| \z-\mu \right\|^2_{\R^{-1}} \right)}{\left\| \z-\mu \right\|^{(n-2)/2}_{\R^{-1}}},
\end{split}
\end{equation}
where $K_n$ is the modified Bessel function of the second kind and order $n$. From \eqref{eqn:laplace-mv} we have that
\[
\tt{E}\left[ \z \right] = \mu, \quad \tt{Cov}\left[ \z \right] = \lambda\, \R.
\]
For large deviations $\left\| \z-\mu \right\|_{\R^{-1}} \to \infty$ the pdf \eqref{eqn:laplace-mv} is approximated by
\begin{equation}
\label{eqn:laplace-mv-approximation}
\mathcal{P}(\z) \approx \left\| \z-\mu \right\|_{\R^{-1}}^{(1-n)/2}\; \exp\left( - \frac{2}{\lambda}\, \left\| \z-\mu \right\|_{\R^{-1}} \right).
\end{equation}
Let $\z$ be the scaled innovation vector \eqref{eqn:define-z}, with $\tt{E}\left[ \z \right] = 0$, and $\tt{Cov}\left[ \z \right] = \mathbf{I}$. The observation part of the variational cost function \eqref{eqn:3dvar-L1} is the log-likelihood, which for large deviations is approximately equal to:
\begin{equation}
\label{eqn:laplace-mv-innovation}
\log\; \mathcal{P}(\z) \approx \frac{n-1}{2}\, \log\, \left\| \z \right\|_2 + \frac{2}{\lambda}\, \left\| \z \right\|_2.
\end{equation}

\end{document}